\documentclass[leqno,10pt]{amsart}
\usepackage{epsfig,latexsym,amssymb,amsmath,graphics,color,subcaption}
\usepackage[foot]{amsaddr}
\pdfoutput=1
\newtheorem{lemma}{Lemma}[section]

\newtheorem{theorem}[lemma]{Theorem}
\newtheorem{proposition}[lemma]{Proposition}
\newtheorem{definition}[lemma]{Definition}
\newtheorem{corollary}[lemma]{Corollary}
\newtheorem{example}[lemma]{Example}
\newtheorem{exercise}[lemma]{Exercise}
\newtheorem{remark}[lemma]{Remark}
\newtheorem{fig}[lemma]{Figure}
\newtheorem{tab}[lemma]{Table}

\newcommand{\bth}{\begin{theorem}}
\newcommand{\ethe}{\end{theorem}}
\newcommand{\bre}{\begin{remark}\rm }
\newcommand{\ere}{\end{remark}}
\newcommand{\ble}{\begin{lemma}}
\newcommand{\ele}{\end{lemma}}
\newcommand{\bde}{\begin{definition}}
\newcommand{\ede}{\end{definition}}
\newcommand{\bco}{\begin{corollary}}
\newcommand{\eco}{\end{corollary}}

\newcommand{\bpr}{\begin{proposition}}
\newcommand{\epr}{\end{proposition}}

\newcommand{\bexer}{\begin{exercise}}
\newcommand{\eexer}{\end{exercise}}

\newcommand{\bexam}{\begin{example}\rm}
\newcommand{\eexam}{\end{example}}

\newcommand{\bfi}{\begin{fig}}
\newcommand{\efi}{\end{fig}}

\newcommand{\btab}{\begin{tab}}
\newcommand{\etab}{\end{tab}}

\def\B_e{B_{\eta}(e)}

\newcommand{\1}{{\bf 1}}

\newcommand{\diag}{\operatorname{diag}}
\newcommand{\bfW}{{\mathbf W}}
\newcommand{\un}{\underline}

\newcommand{\twonorm}[1]{\|#1\|_2}

\newcommand{\frobnorm}[1]{\|#1\|_F}

\newcommand{\ov}{\overline}
\newcommand{\wt}{\widetilde}

\definecolor{darkblue}{rgb}{.1, 0.1,.8}
\definecolor{darkgreen}{rgb}{0,0.8,0.2}
\definecolor{darkred}{rgb}{.8, .1,.1}

\newcommand{\sv}{stochastic volatility}

\newcommand{\slln}{strong law of large numbers}
\newcommand{\clt}{central limit theorem}

\newcommand{\asy}{asymptotic}

\newcommand{\ts}{time series}

\newcommand{\bfV}{{\bf V}}

\newcommand{\pp}{point process}

\newcommand{\lhs}{left-hand side}
\newcommand{\fidi}{finite-dimensional distribution}
\newcommand{\rv}{random variable}

\newcommand{\var}{{\rm var}}

\newcommand{\cov}{{\rm cov}}

\newcommand{\bfTh}{\mbox{\boldmath$\Theta$}}

\newcommand{\rhs}{right-hand side}

\newcommand{\beao}{\begin{eqnarray*}}
\newcommand{\eeao}{\end{eqnarray*}\noindent}

\newcommand{\beam}{\begin{eqnarray}}
\newcommand{\eeam}{\end{eqnarray}\noindent}

\newcommand{\beqq}{\begin{equation}}
\newcommand{\eeqq}{\end{equation}\noindent}

\newcommand{\bce}{\begin{center}}
\newcommand{\ece}{\end{center}}

\newcommand{\barr}{\begin{array}}
\newcommand{\earr}{\end{array}}

\newcommand{\stp}{\stackrel{\P}{\rightarrow}}
\newcommand{\std}{\stackrel{d}{\rightarrow}}

\newcommand{\stv}{\stackrel{v}{\rightarrow}}
\newcommand{\stw}{\stackrel{w}{\rightarrow}}

\newcommand{\eqd}{\stackrel{d}{=}}

\newcommand{\vague}{\stackrel{\lower0.2ex\hbox{$\scriptscriptstyle
                    \it{v} $}}{\rightarrow}}
\newcommand{\weak}{\stackrel{\lower0.2ex\hbox{$\scriptscriptstyle
                    \it{w} $}}{\rightarrow}}
\newcommand{\what}{\stackrel{\lower0.2ex\hbox{$\scriptscriptstyle
                    \it{\hat{w}} $}}{\rightarrow}}

\newcommand{\bdis}{\begin{displaymath}}
\newcommand{\edis}{\end{displaymath}\noindent}

\newcommand{\nto}{n\to\infty}

\newcommand{\xto}{x\to\infty}

\newcommand{\wh}{\widehat}
\newcommand{\vep}{\varepsilon}

\newcommand{\la}{\lambda}

\newcommand{\regvary}{regularly varying}
\newcommand{\slvary}{slowly varying}
\newcommand{\regvar}{regular variation}

\newcommand{\bbr}{{\mathbb R}}

\newcommand{\bbz}{{\mathbb Z}}

\renewcommand{\P}{{\mathbb P}}
\newcommand{\E}{{\mathbb E}}

\newcommand{\con}{convergence}

\newcommand{\evt}{extreme value theory}

\newcommand{\st}{such that}
\newcommand{\fif}{if and only if}

\newcommand{\chf}{characteristic function}
\newcommand{\fct}{function}

\newcommand{\ds}{distribution}

\newcommand{\rep}{representation}

\newcommand{\seq}{sequence}

\newcommand{\pro}{probabilit}

\newcommand{\ms}{measure}

\newcommand{\bfx}{{\bf x}}
\newcommand{\bfX}{{\bf X}}

\newcommand{\bfY}{{\mathbf Y}}

\newcommand{\bfA}{{\bf A}}

\newcommand{\bfZ}{{\bf Z}}

\newcommand{\bfS}{{\bf S}}

\newcommand{\bfe}{{\bf e}}

\newcommand{\bft}{{\bf t}}
\newcommand{\bfI}{{\bf I}}
\newcommand{\bfs}{{\bf s}}

\newcommand{\ex}{{\rm e}}
\begin{document}
\today
\title[Eigenvalues of the sample covariance matrix of a   
stochastic volatility model]{The eigenvalues of the sample covariance matrix of a multivariate heavy-tailed 
stochastic volatility model}
%\author[Johannes Heiny]{Johannes Heiny}
\author{Anja Jan\ss en$^{1,\ast}$}
\author{Thomas Mikosch$^1$} 
\author{Mohsen Rezapour$^2$} 
\author{Xiaolei Xie$^1$}
\address{$^{1}$ Department  of Mathematical Sciences,
University of Copenhagen,
Universitetsparken 5,
DK-2100 Copenhagen,
Denmark}
\address{$^2$
Department of Statistics, Faculty of Mathematics and Computer, Shahid
Bahonar University of Kerman, Kerman, Iran}
\thanks{$^\ast$ Corresponding author}
\email{anja@math.ku.dk}
\email{mikosch@math.ku.dk\, (URL: www.math.ku.dk/$\sim$mikosch)}
\email{mohsenrzp@gmail.com}
\email{xie.xiaolei@gmail.com}

\begin{abstract}
We consider a multivariate heavy-tailed stochastic volatility model and analyze the large-sample behavior of its sample covariance matrix. We study the limiting behavior of its entries in the infinite-variance case and derive results for the ordered eigenvalues and corresponding eigenvectors. Essentially, we consider two different cases where the tail behavior either stems from the iid\ innovations of the process or from its volatility sequence. In both cases, we make use of a large deviations technique for regularly varying time series to derive multivariate $\alpha$-stable limit distributions of the sample covariance matrix. While we show that in the case of heavy-tailed innovations the limiting behavior resembles that of completely independent observations, we also derive that in the case of a heavy-tailed volatility sequence the possible limiting behavior is more diverse, i.e.\ allowing for dependencies in the limiting distributions which are determined by the structure of the underlying volatility sequence. 
\end{abstract}


\maketitle

\section{Introduction}\label{sec:intro}
\subsection{Background and Motivation}\label{subsec:motivation}

The study of sample covariance matrices is fundamental for the analysis of dependence in multivariate time series. Besides from providing estimators for variances and covariances of the observations (in case of their existence), the sample covariance matrices are a starting point for dimension reduction methods like principal component analysis. Accordingly, the special structure of sample covariance matrices and their largest eigenvalues has been intensively studied in random matrix theory, starting with iid\ Gaussian observations and more recently extending results to arbitrary distributions which satisfy some moment assumptions like in the four moment theorem of Tao and Vu \cite{tao:vu:2012}. 

However, with respect to the analysis of financial time series, such a moment assumption is often not suitable. Instead, in this work, we will analyze the large sample behavior of sample covariance matrices under the assumption that the marginal distributions of our observations are regularly varying with index $\alpha<4$ which implies that fourth moments do not exist. In this case, we would expect the largest eigenvalues of the sample covariance matrix to inherit 
heavy-tailed behavior as well; see for example Ben Arous and 
Guionnet \cite{benarous:guionnet:2007}, 
Auffinger et al. \cite{auffinger:arous:peche:2009}, Soshnikov \cite{soshnikov:2004,soshnikov:2006}, Davis et al. \cite{davis:et:al:2016}, Heiny and Mikosch \cite{heiny:mikosch:2016}
for the case of iid entries. Furthermore, in the context of financial time series we have to allow for dependencies both over time and between different components and indeed it is the very aim of the analysis to discover and test for these dependencies from the resulting sample covariance matrix as has for example been done in 
Plerou et al. \cite{plerou:et:al:2002} and Davis et al. \cite{davis:pfaffel:stelzer:2014,davis:mikosch:pfaffel:2016}. 
The detection of dependencies among assets also plays a crucial role in portfolio optimization based on multi-factor prizing models, where principal component analysis is one way to derive the main driving factors of a portfolio; cf.\ Campbell et al. \cite{campbell:lo:macKinlay:1997} and recent 
work by Lam and Yao \cite{lam:yao:2012}.
\par
The literature on the asymptotic behavior of sample covariance matrices derived from dependent heavy-tailed data is, however, relatively sparse up till now. Starting with the analysis of the sample autocorrelation of univariate linear heavy-tailed time series in Davis and Resnick \cite{davis:resnick:1985, davis:resnick:1986}, 
the theory has recently been extended to multivariate heavy-tailed time series with linear structure in 
Davis et al. \cite{davis:pfaffel:stelzer:2014,davis:mikosch:pfaffel:2016}, cf.\ also the recent 
survey article by Davis et al. \cite{davis:et:al:2016}. But most of the standard models for financial time series
such as GARCH and \sv\ models
are non-linear. In this paper we will therefore focus on a class of multivariate stochastic volatility models of the form 
\beam\label{Eq:intro:sv}
X_{it}=\sigma_{it}\,Z_{it}\,,\qquad t\in\bbz\,, \quad 1\leq i \leq p,
\eeam
where $(Z_{it})$ is an iid random field independent of a strictly stationary ergodic field $(\sigma_{it})$ of non-negative \rv s; 
see Section \ref{sec:model} for further details.  Stochastic volatility models have been studied in detail in the financial \ts\ literature; see for example
Andersen et al. \cite{andersen:davis:kreiss:mikosch:2009}, Part II. They are among the simplest models allowing for
conditional heteroscedasticity of a time series. In view of independence between the $Z$- and $\sigma$-fields 
dependence conditions on $(X_{it})$ are imposed only via the \sv\ $(\sigma_{it})$. Often it is assumed that $(\log \sigma_{it})$ has a linear 
structure, most often Gaussian.
\par
In this paper we are interested in the case when the marginal and \fidi s of $(X_{it})$ have power-law tails.
Due to independence between  $(\sigma_{it})$ and $(Z_{it})$ heavy tails of $(X_{it})$ can be due to the $Z$- or the $\sigma$-field.
Here we will consider two cases: (1) the tails of $Z$ dominate the right tail of $\sigma$ and (2) the right tail of 
$\sigma$ dominates  the tail of $Z$. The third case when both $\sigma$ and $Z$ have heavy tails and are tail-equivalent
will not be considered in this paper. Case (1) is typically more simple to handle; see Davis and Mikosch 
\cite{davis:mikosch:1999,davis:mikosch:2001,davis:mikosch:2009} for \evt , \pp\ \con\ and central limit theory 
with infinite variance stable limits. Case (2) is more subtle as regards the tails of the \fidi s. The literature on \sv\ models with a heavy-tailed volatility sequence is so far sparse but the interest in these models has been growing recently; see  
Mikosch and Rezapour \cite{mikosch:rezapour:2013}, Kulik and Soulier \cite{kulik:soulier:2015} and Jan\ss en and Drees \cite{janssen:drees:2016}. In particular, it has been shown that these models offer a lot of flexibility with regard to the extremal dependence structure of the time series, ranging from asymptotic dependence of consecutive observations (cf.\ \cite{mikosch:rezapour:2013}) to asymptotic independence of varying degrees (cf.\ \cite{kulik:soulier:2015} and \cite{janssen:drees:2016}). 
\subsection{Aims, main results and structure} 
After introducing the general model in Section \ref{sec:model} we first deal with the case of heavy-tailed innovations and a light-tailed volatility sequence in Section \ref{sec:case1}. The first step in our analysis is to describe the extremal structure of the corresponding process by deriving its so-called tail process; see  
Section~\ref{subsec:rvts} and Proposition~\ref{prop:regvarsv}. This allows one to apply an infinite variance 
stable \clt\  from Mikosch and Wintenberger \cite{mikosch:wintenberger:2016} (see Appendix \ref{App:A}) to derive the joint limiting behavior of the entries of the sample covariance matrix of this model. This leads to the main results 
in the first case: Theorems~\ref{the:1} and \ref{lem:1}. They say, roughly speaking, that all values on the off-diagonals of the sample covariance matrix are negligible compared to the values on the diagonals.
Furthermore, the values on the diagonal converge, under suitable normalization, to independent $\alpha$-stable random variables, so the limiting behavior of this class of 
stochastic volatility  models is quite similar to the case of iid\ heavy-tailed random variables. This fairly tractable structure allows us also to derive explicit results about the asymptotic behavior of the ordered eigenvalues and corresponding eigenvectors which can be found in Sections \ref{subsec:eigen} and \ref{subsec:appl}. In particular, we will see that in this model 
the eigenvectors are basically the unit canonical basis vectors which describe a very weak form of extremal dependence. With a view towards portfolio analysis, our assumptions imply that large movements of the market are mainly driven by one single asset, where each asset is equally likely to be this extreme driving force. 

In the second case of a heavy-tailed volatility sequence combined with light-tailed innovations, which we analyze in Section \ref{Sec:case2}, we see that the range of possible limiting behaviors of the entries of the sample covariance matrix is more diverse and depends on the specific structure of the underlying volatility process. We make the common assumption that our volatility process is log-linear, where we distinguish between two different cases for the corresponding innovation distribution of this process. Again, for both cases, we first derive the specific form of the corresponding tail process (see Proposition~\ref{prop:rv:2}) 
which then allows us to derive the limiting behavior of the sample covariance matrix entries, leading to the main results in the second case: Theorems~\ref{thm:1} and \ref{thm:2}. We show that the sample covariance matrix can feature non-negligible off-diagonal components, therefore clearly distinguishing from the iid\ case, if we assume that the innovations of the log-linear volatility process are convolution equivalent. We discuss concrete examples for both model specifications and the corresponding implications for the asymptotic behavior of ordered eigenvalues and corresponding eigenvectors at the end of Section \ref{Sec:case2}. 

Section \ref{sec:simulation} contains a small simulation study which illustrates our results for both cases and also includes a 
real-life data example for comparison. From the foreign exchange rate data that we use, it is notable that the corresponding sample covariance matrix features a relatively large gap between the largest and the second largest eigenvalue and that the eigenvector corresponding to the largest eigenvalue is fairly spread out, i.e.,\ all its components are of a similar order of magnitude. This implies that the model discussed in Section~\ref{sec:case1} may not be that suitable to catch the extremal dependence of this data, and that there is not one single component that is most affected by extreme movements but instead all assets are 
affected in a similar way. We perform simulations for three different specifications of models 
from Sections~\ref{sec:case1} and \ref{Sec:case2}. They illustrate that the models analyzed 
in Section~\ref{Sec:case2} are capable of exhibiting more diverse asymptotic behaviors 
of the sample covariance matrix and in particular non-localized dominant eigenvectors.

Some useful results for the (joint) tail and extremal behavior of random products are gathered in Appendix~\ref{App:B}. 
These results may be of independent interest when studying the extremes of multivariate \sv\ models with possibly distinct tail
indices. We mention in passing that there is great interest in non-linear models for log-returns of speculative prices when
the number of assets $p$ increases with the sample size $n$. We understand our analysis as a first step in this direction.

\section{The model}\label{sec:model}\setcounter{equation}{0}
We consider a \sv\ model
\beam\label{eq:1a}
X_{it}=\sigma_{it}\,Z_{it}\,,\qquad i,t\in\bbz\,,
\eeam
where $(Z_{it})$ is an iid field independent of a strictly stationary ergodic field $(\sigma_{it})$ of non-negative \rv s.
We write $Z$, $\sigma$, $X$ for generic elements of the $Z$-, $\sigma$- and $X$-fields \st\ $\sigma$ and $Z$ are independent.
A special case appears when $\sigma>0$ is a constant: then $(X_{it})$ constitutes an iid field. 

For the \sv\ model as in \eqref{Eq:intro:sv} we construct the 
multivariate \ts\ 
\beam\label{eq:multvts}
 \bfX_t= (X_{1t},\ldots,X_{pt})', \;\;\; t \in \mathbb{Z},
\eeam 
for a given dimension $p\ge 1$. For $n\geq 1$ we write $\bfX^n={\rm vec}\big((\bfX_t)_{t=1, \ldots, n}\big) \in \mathbb{R}^{p \times n}$ and consider the non-normalized 
sample covariance matrix
\beam\label{eq:esses}
\bfX^n(\bfX^n)'= (S_{ij})_{i,j=1,\ldots,p}\,,\qquad S_{ij}= \sum_{t=1}^n X_{it}X_{jt}\,,\qquad S_i=S_{ii}\,.
\eeam

\subsection{Case (1): $Z$ dominates the tail}
We assume that $Z$ is \regvary\ with index $\alpha>0$, i.e.,
\beam\label{eq:10b}
\P(Z>x)\sim p_+\,\dfrac{L(x)}{x^\alpha}\qquad \mbox{and}\qquad \P(Z<-x)\sim p_-\,\dfrac{L(x)}{x^\alpha}\,,\qquad \xto\,,
\eeam
where $p_+$ and $p_-$ are non-negative numbers with $p_++p_-=1$ and $L$ is a \slvary\ \fct . 
If we assume $\E[\sigma^{\alpha+\delta}]<\infty$ for some $\delta>0$ then, in view of a result 
by Breiman \cite{breiman:1965} (see also Lemma~\ref{lem:product}), it follows that
\beam\label{eq:br1}
\P(X>x)\sim \E[\sigma^\alpha]\,\P(Z>x)\quad \mbox{and}\quad \P(X<-x)\sim \E[\sigma^\alpha]\,\P(Z<-x)\,,\qquad \xto\,, 
\eeam
i.e., $X$ is \regvary\ with index $\alpha$. Moreover, we know from a result by Embrechts and Goldie \cite{embrechts:goldie:1980}
that for independent copies $Z_1$ and $Z_2$ of $Z$, $Z_1Z_2$ is again \regvary\ with index $\alpha$; cf. Lemma~\ref{lem:product}. Therefore, using again
Breiman's result under the 
condition that $\E[(\sigma_{i0}\sigma_{j0})^{\alpha+\delta}\1(i\ne j)+\sigma_{i0}^{\alpha+\delta}]<\infty$ for some $\delta>0$, we have 
\beam\label{eq:br2}
\P(\pm X_{it}\, X_{jt}>x) \sim \left\{\begin{array}{ll}\E[(\sigma_{it}\,\sigma_{jt})^{\alpha}]\,\P(\pm Z_i\,Z_j>x)& i\ne j\,,\\[2mm]
 \E[\sigma^{\alpha}]\,\P(Z^2>x)& i=j\,,
\end{array}\qquad \xto\,.
\right.
\eeam 
\subsection{Case (2): $\sigma$ dominates the tail}\label{subsec:Case:2}
We assume that $\sigma\ge 0$ is \regvary\ with some index $\alpha>0$: for some \slvary\ \fct\ $\ell$,
\beao
\P(\sigma>x)= x^{-\alpha}\,\ell(x)\,,
\eeao 
and $\E [|Z|^{\alpha+\delta}]<\infty$ for some $\delta>0$. Now the Breiman result yields
\beao
\P(X>x) \sim \E [Z_+^\alpha]\,\P(\sigma>x)\qquad \mbox{and} \qquad \P(X<-x) \sim \E[Z_-^\alpha]\,\P(\sigma>x)\,,\qquad \xto\,.
\eeao
Since we are also interested in the tail behavior of the products $X_{it}X_{jt}$ we need to
be more precise about the joint \ds\ of the \seq s $(\sigma_{it})$. We assume
\beam\label{eq:sv2model}
\sigma_{it}= \exp\Big(\sum_{k,l=-\infty}^\infty \psi_{kl} \,\eta_{i-k,t-l}\Big)\,,\qquad i,t\in\bbz\,,
\eeam
where $(\psi_{kl})$ is a field of non-negative numbers (at least one of them being positive) \st\ (without loss of generality)
$\max_{kl}\psi_{kl}=1$
and 
$(\eta_{it})$ is an iid random field  \st\ a generic element $\eta$ satisfies
\beam\label{Eq:tail:eta}
\P\big(\ex^{\eta}>x)= x^{-\alpha}\,L(x)\,, 
\eeam
for some $\alpha>0$ and a \slvary\ \fct\ $L$. 
We also assume  
$\sum_{k,l}\psi_{kl}<\infty$ to ensure absolute summability of 
$\log \sigma_{it}$. 
A distribution of $\eta$ that fits into this scheme is for example the exponential distribution; cf.\ also Rootz\'{e}n \cite{rootzen:1986} for further examples and extreme value theory for linear processes of the 
form $\sum_{l=-\infty}^\infty \psi_{l} \,\eta_{t-l}$. 
\subsection{Regularly varying \seq s}\label{subsec:rvts}
In Sections~\ref{subsec:regvar} and \ref{subsec:regvar2} 
we will elaborate on the joint tail behavior of the \seq s $(\sigma_{it})$, $(X_{it})$, $(\sigma_{it}\sigma_{jt})$, 
and $(X_{it}X_{jt})$. We will show that, under suitable conditions,  these \seq s are \regvary\ with  
positive indices. 
\par
The notion of a {\em univariate \regvary\ \seq } was introduced by Davis and Hsing \cite{davis:hsing:1995}. 
Its extension to the multivariate case does not represent difficulties; see Davis and Mikosch \cite{davis:mikosch:2009b}.
An $\bbr^d$-valued strictly stationary \seq\ $(\bfY_t)$ is {\em \regvary\ with index $\gamma>0$} if each of the vectors
$(\bfY_t)_{t=0,\ldots,h}$, $h\ge 0$, is \regvary\ with index $\gamma$, i.e., there exist non-null Radon \ms s $\mu_h$ on 
$[-\infty,\infty]^{d(h+1)}\backslash\{\bf0\}$ which are homogeneous of order $-\gamma$ \st
\beam \label{eq:multvtsvague}
\dfrac{\P(x^{-1} (\bfY_t)_{t=0,\ldots,h}\in \cdot)}{\P(\|\bfY_0\|>x)} \stv \mu_h(\cdot)\,.
\eeam
Here $\stv$ denotes vague \con\ on the Borel $\sigma$-field of $[-\infty,\infty]^{d(h+1)}\backslash\{\bf0\}$ and $\|\cdot\|$ denotes any given norm; see Resnick's books \cite{resnick:1987,resnick:2007} as general references to multivariate \regvar .
\par
Following Basrak and Segers \cite{basrak:segers:2009}, 
an $\bbr^d$-valued strictly stationary \seq\ $(\bfY_t)$ is \regvary\ with index $\gamma>0$ \fif\ there exists a \seq\ 
of $\bbr^d$-valued random vectors $(\bfTh_h)$ independent of a Pareto($\gamma$) \rv\ 
$Y$, i.e., $\P(Y>x)=x^{-\gamma}$, $x>1$,  \st\ for any $k\ge 0$, 
\beam\label{eq:basseg}
\P( x^{-1} (\bfY_0,\ldots,\bfY_k)\in \cdot \mid \|\bfY_0\|>x)\stw \P\big(Y\,(\bfTh_0,\ldots,\bfTh_k)\in\cdot\big)\,,\qquad
\xto\,.
\eeam
We call $(\bfTh_h)$ the {\em spectral tail process} of $(\bfY_t)$ and $(Y \bfTh_h)$ the {\em tail process}. We will use both defining properties (i.e.,\ \eqref{eq:multvtsvague} and \eqref{eq:basseg}) of a \regvary\ \seq .

\section{Case (1): $Z$ dominates the tail}\label{sec:case1}\setcounter{equation}{0}
\subsection{Regular variation of the \sv\ model and its product processes}\label{subsec:regvar}
\bpr\label{prop:regvarsv}
We assume the \sv\ model \eqref{eq:1a} and that $Z$
is \regvary\ with index $\alpha>0$ in the sense of \eqref{eq:10b}.% and  $\E[\sigma^{\alpha+\vep}]<\infty$ for some $\vep>0$.
\begin{enumerate}
\item
If $\E[\sigma^{\alpha+\vep}]<\infty$ for some $\vep>0$
the \seq\ $(X_{it})_{t \in\bbz}$ is \regvary\ with index $\alpha$ and the corresponding spectral tail process $(\Theta^i_h)_{h\ge 1}$ vanishes.
\item
For any $i\ne j$, if $\E[(\sigma_{i0}\sigma_{j0})^{\alpha+\vep}]<\infty$ for some $\vep>0$ then the \seq\ $(X_{it}X_{jt})$ is \regvary\ with index $\alpha$ and the corresponding spectral tail process 
$(\Theta_h^{ij})_{h\ge 1}$ vanishes.
\end{enumerate}
\epr
\bre
If  $\E[(\sigma_{ik}\sigma_{jl})^{\alpha+\vep_{ik,jl}}]<\infty$ for some $\vep_{ik,jl}>0$ and any $(i,k)\ne (j,l)$ it is also possible to show 
the joint \regvar\ of the processes $(X_{it}X_{jt})$, $i\ne j$, with index $\alpha$. The description of the corresponding spectral tail process
is slightly tedious. It is not needed for the purposes of this paper and therefore omitted.
\ere
\begin{proof} Regular variation of the marginal \ds s of  $(X_{it})$ and $(X_{it}X_{jt})$ follows from Breiman's result; see
\eqref{eq:br1} and \eqref{eq:br2}. As regards the \regvar\ of the \fidi s of $(X_{it})$, we have for $h\ge 1$,
\beao
\P(|X_{ih}|>x\mid |X_{i0}|>x)&=&\dfrac{\P(\min(|X_{i0}|,|X_{ih}|)>x)}{\P(|X_{i0}|>x)}\\
&\le &\dfrac{\P(\max(\sigma_{i0},\sigma_{ih}) \min(|Z_{i0}|,|Z_{ih}|)>x)}{\P(|X_{i0}|>x)}\to 0\,,\qquad \xto\,.
\eeao
In the last step we used Markov's inequality together with the moment condition $\E[\sigma^{\alpha+\vep}]<\infty$ and the fact that 
$\min(|Z_{i0}|,|Z_{ih}|)$ is \regvary\ with index $2\alpha$.
This means that $\Theta^i_h=0$ for $h\ge 1$. 
\par
Similarly, for $i\ne j$, $h\ge 1$,
\beao
\P(|X_{ih}X_{jh}|>x\mid |X_{i0}X_{j0}|>x)&\le &\dfrac{\P(\max(\sigma_{i0}\sigma_{j0},\sigma_{ih}\sigma_{jh}) \min(|Z_{i0}Z_{j0}|,
|Z_{ih}Z_{jh}|)>x)}{\P(|X_{i0}X_{j0}|>x)}\to 0\,.
\eeao
In the last step we again used Markov's inequality, the fact that $Z_{i0}Z_{j0}$ is \regvary\ with index $\alpha$
(see  Embrechts and Goldie \cite{embrechts:goldie:1980}; cf.\ Lemma~\ref{lem:product}(1) below), hence $\min(|Z_{i0}Z_{j0}|,|Z_{ih}Z_{jh}|)$ is \regvary\ with index~$2\alpha$, and the moment condition $\E [(\sigma_{i0}\sigma_{j0})^{\alpha+\vep}]<\infty$.
Hence $\Theta_h^{ij}=0$ for $i\ne j$, $h\ge 1$.
%The joint \regvar\ of the processes $(X_{it}X_{jt})$, $i\ne j$, follows in a similar. 
%The same arguments as above show that for $i\ne j$, $k\le l$,
%\beao
%\P(|X_{ih}X_{jh}|>x\mid |X_{k0}X_{l0}|>x)\to 0\,,\qquad \xto\,.
%\eeao
%Some care has to be taken if $k=i$. Then
%\beao
%\P( X_{ih}X_{jh}>x\mid |X_{i0}X_{l0}|>x)\to 0
%\eeao
\end{proof}
%as well as sample correlation matrix
%\beao
%\wt\bfX\wt \bfX'= \big(\sum_{t=1}^n\wt X_{it}\wt X_{jt}\big)_{i,j=1,\ldots,p}\,,
%\eeao
%where 
%\beao
%\wt X_{it}= \dfrac{X_{it}- \ov X_i}{\wh \sigma_i}\,,\qquad \ov X_i= \dfrac 1n \sum_{t=1}^n X_{it}\,,\qquad 
%\wh \sigma_i^2= \dfrac 1n \sum_{t=1}^n (X_{it}-\ov X_i)^2\,. 
%\eeao
%We notice that ${\rm diag}(\wt\bfX\wt\bfX')=\bfI_p$, the $p\times p$ identity matrix.
\subsection{Infinite variance stable limit theory for the \sv\ model and its product processes}\label{subsec:limit}
\bth\label{the:1}
Consider the \sv\ model \eqref{eq:1a} and assume the following conditions:
\begin{enumerate}
\item
$Z$ is \regvary\ with index $\alpha\in (0,4) \setminus \{2\}$.
\item
$\big((\sigma_{it})_{t=1,2,\ldots}\big)_{i=1,\ldots,p}$  is strongly mixing with rate \fct\ $(\alpha_h)$ \st\ for some $\delta>0$,
\beam\label{eq:mixi}
\sum_{h=0}^\infty \alpha_h^{\delta/(2+\delta)}<\infty\,.
\eeam
\item
The moment condition
\beam\label{eq:moment}
\E[\sigma^{2\max(2+\delta,\alpha+\epsilon)} ]<\infty
\eeam
holds for the same $\delta>0$ as in \eqref{eq:mixi} and some $\epsilon>0$. 
\end{enumerate}
Then 
\beam\label{eq:feller}
a_n^{-2}\big(S_{1}-c_n,\ldots,S_{p}-c_n\big) \std (\xi_{1,\alpha/2},\ldots,\xi_{p,\alpha/2})\,,
\eeam
where $(\xi_{i,\alpha/2})$ are iid $\alpha/2$-stable \rv s which are totally skewed to the right,
\beam\label{eq:6}
c_n=\left\{\begin{array}{ll}
0& \alpha \in (0,2)\,,\\
n \,\E[X^2] &\alpha\in (2,4)\,,
\end{array}
\right.
\eeam
and $(a_n)$ satisfies $n\,\P(|X|>a_n)\to 1$ as $\nto$.
\ethe
\bre\label{rem:feller}
From classical limit theory (see Feller \cite{feller:1971}, Petrov \cite{petrov:1995}) we know 
that \eqref{eq:feller} holds for an iid random field $(X_{it})$ with \regvary\ $X$ with index $\alpha\in (0,4)$. In the case $\alpha=2$ one needs the special centering $c_n=n\,\E [X^2 \1(|X|\le a_n)]$ which often leads to some additional
technical difficulties. For this reason we typically exclude this
case in the sequel.
\ere
\bre\label{rem:mix}
It follows from standard theory that $\alpha$-mixing of $(\sigma_{it})$ with rate \fct\ $(\alpha_h)$ implies $\alpha$-mixing of 
$(X_{it})$ with rate \fct\ $(4\alpha_h)$; see Davis and Mikosch \cite{davis:mikosch:2009}.
\ere
\begin{proof} 
Recall the definition of $(\bfX_t)$ from \eqref{eq:multvts}. We will verify the conditions of Theorem~\ref{thm:mikwin}
for $\bfX_t^2=(X_{it}^2)_{i=1,\ldots,p}$, $t=0,1,2,\ldots$.\\[2mm]
(1) We start by verifying the \regvar\ condition for $(\bfX_t)$;  see  \eqref{eq:basseg}.
We will determine the \seq\ $(\bfTh_h)$ corresponding to $(\bfX_t)$.
We have for $t\ge 1$, with the max-norm $\|\cdot\|$,
\beao
\P\big(\|\bfX_t\|> x \mid \|\bfX_0\|>x\big)&\le &
\dfrac{\P\big( \|\bfX_t\|> x\,, \cup_{i=1}^p\{|X_{i0}|>x\}\big)}{\P(\|\bfX_0\|>x)}\\
&\le &\sum_{i=1}^p \dfrac{\P\big(  \|\bfX_t\|> x\,, |X_{i0}|>x\big)}{\P(\|\bfX_0\|>x)}\\
&\le &\sum_{i=1}^p \sum_{j=1}^p\dfrac{\P\big( |X_{jt}|> x\,, |X_{i0}|>x\big)}{\P(|X|>x)}\\
&\le &\sum_{i=1}^p \sum_{j=1}^p\dfrac{\P\big( \max(\sigma_{jt},\sigma_{i0})\min(|Z_{jt}|,|Z_{i0}|)> x\big)}{\P(\sigma|Z|>x)}\,.
\eeao
We observe that by Breiman's result and in view of the moment condition \eqref{eq:moment}, for $t\ge 1$ and some positive constant $c$,
\beao
\dfrac{\P\big( \max(\sigma_{jt},\sigma_{i0})\min(|Z_{jt}|,|Z_{i0}|)> x\big)}{\P(\sigma|Z|>x)} \sim c\, \dfrac{\P(\min(|Z_{jt}|,|Z_{i0}|)> x)}{\P(|Z|>x)}\,,
\eeao
and the \rhs\ converges to zero as $\xto$.
We conclude that $\bfTh_h=\bf0$ for $h\ge 1$. We also have for $i\ne j$,
\beao
\dfrac{\P(|X_{i0}|>x\,,|X_{j0}|>x)}{\P(|X|>x)} \leq \dfrac{\P\big( \max(\sigma_{i0},\sigma_{j0})\min(|Z_{i0}|,|Z_{j0}|)> x\big)}{\P(\sigma|Z|>x)} \to 0\,,\qquad \xto \,.
\eeao
Then, in a similar way, one can show
\beam\label{eq:Theta}
\P(\bfX_0/\|\bfX_{0}\|\in \cdot \mid \|\bfX_0\|>x)
&\stw & \P(\bfTh_0 \in \cdot) =\dfrac  1p \sum_{i=1}^p \big(p_+\vep_{\bfe_i}(\cdot)+p_-\vep_{-\bfe_i}(\cdot)\big)\,.
\eeam
where $\bfe_i$ are the canonical basis vectors in $\bbr^p$, $\vep_\bfx$ is Dirac \ms\ at $\bfx$ and $p_\pm$ are the tail
balance factors in \eqref{eq:10b}. 
\par
We conclude that the spectral tail process $(\bfTh_h^{(2)})$ of $(\bfX_t^2)$ is given by 
$\bfTh_h^{(2)}=\bf0$ for $h\ge 1$ and from \eqref{eq:Theta} we also have
\beam\label{Eq:Theta_0:measure}
\P(\bfTh_0^{(2)}\in\cdot) =\dfrac  1p \sum_{i=1}^p \vep_{\bfe_i}(\cdot)\,.
\eeam
In particular, the condition $\sum_{i=1}^\infty  \E [\|\bfTh_i^{(2)}\|]<\infty$  in Theorem~\ref{thm:mikwin}(4)
is trivially satisfied.\\[2mm]
(2) Next we want to prove the mixing condition~\eqref{eq:chfa} for the \seq\ $(\bfX_t^2)$.  
We start by observing that there are integer \seq s  $(l_n)$ and $(m_n)$ \st\ $k_n\,\alpha_{l_n}\to 0$, $l_n=o(m_n)$ and $m_n=o(n)$. Then 
we  also have for any $\gamma>0$,
\beam\label{eq:neg}
k_n\,\P\big( \sum_{t=1}^{l_n} \bfX_t^2 \1(\|\bfX_t\|>\vep a_n)>\gamma a_n^2 \big)\le k_n\,l_n\,\P(\|\bfX_t\|>\vep a_n)\le c\, l_n/m_n=o(1)\,.\eeam
Relation \eqref{eq:chfa} turns into
\beao
\E \ex^{i \bfs 'a_n^{-2} \sum_{t=1}^n \bfX_t^2 \1(\|\bfX_t\|>\vep a_n)}- \big(\E \ex^{i\bfs' a_n^{-2}\sum_{t=1}^{m_n} \bfX_t^2 \1(\|\bfX_t\|>\vep a_n)}\big)^{k_n}\to 0\,,\qquad \bfs\in\bbr^p\,.
\eeao
In view of \eqref{eq:neg} it is not difficult to see that we can replace the sum in the former \chf\ by the sum over the index set
$J_n=\{1,\ldots,m_n-l_n,m_n+1,\ldots,2m_n-l_n,\ldots,\}\subset \{1,\ldots,n\}$ and in the latter \chf\ by the sum over the index set $\{1,\ldots,m_n-l_n\}$. 
Without loss of generality we may assume that $n/m_n$ is an integer. Thus it remains to show that the following difference 
converges to zero for every $\bfs\in\bbr^p$:
\beao\lefteqn{\Big|
\E \big[\ex^{i \bfs 'a_n^{-2} \sum_{t\in J_n} \bfX_t^2 \1(\|\bfX_t\|>\vep a_n)}\big]- 
\Big(\E\big[ \ex^{i\bfs' a_n^{-2}\sum_{t=1}^{m_n-l_n} \bfX_t^2 \1(\|\bfX_t\|>\vep a_n)}\big]\Big)^{k_n}\Big|}\\&=&\Big| 
\sum_{v=1}^{k_n} \E\Big[
\prod_{j=1}^{v-1} \ex^{i \bfs 'a_n^{-2} \sum_{t=(j-1)m_n+1}^{jm_n-l_n} \bfX_t^2 
\1(\|\bfX_t\|>\vep a_n)}\\&&\times  \big(\ex^{i \bfs 'a_n^{-2} \sum_{t=(v-1)m_n+1}^{vm_n-l_n} \bfX_t^2 
\1(\|\bfX_t\|>\vep a_n)} -\E \big[\ex^{i \bfs 'a_n^{-2} \sum_{t=(v-1)m_n+1}^{vm_n-l_n} \bfX_t^2 
\1(\|\bfX_t\|>\vep a_n)}\big] \big)\Big]\\&&\times  \prod_{j=v+1}^{k_n} \E \big[\ex^{i \bfs 'a_n^{-2} \sum_{t=(j-1)m_n+1}^{jm_n-l_n} \bfX_t^2 
\1(\|\bfX_t\|>\vep a_n)}\big]\Big|\,.  
\eeao 
In view of a standard inequality for covariances of strongly mixing \seq s of bounded \rv s (see Doukhan \cite{doukhan:1994}, p.~3) the 
\rhs\ is bounded by $c\,k_n\alpha_{l_n}$ which converges to zero by construction. Here and in what follows, $c$ stands for any 
positive constant whose value is not of interest. Its value may change from line to line.
This finishes the proof of the mixing condition.\\[2mm]
(3)~Next we check the anti-clustering condition \eqref{eq:acl} for $(\bfX_t)$ with normalization $(a_n)$, implying the corresponding
condition for $(\bfX_t^2)$ with normalization $(a_n^2)$.
By similar methods as for part (1) of the proof, assuming that $\|\cdot\|$ is the max-norm, we have 
\beao\lefteqn{
\P\big(\max_{t=l,\ldots,m_n} \|\bfX_t\|>\gamma a_n\mid \|\bfX_0\|>\gamma a_n\big)}\\&\le &
\sum_{t=l}^{m_n} \P\big( \|\bfX_t\|>\gamma a_n\mid \|\bfX_0\|>\gamma a_n\big)\\
&\le &c\,\sum_{t=l}^{m_n} \sum_{i=1}^p\sum_{j=1}^p\dfrac{\P\big( |X_{it}|>\gamma a_n\,,|X_{j0}|>\gamma a_n\big)}{\P(|Z|>\gamma a_n)}\\
&\le &c\,\sum_{t=l}^{m_n} \sum_{i=1}^p\sum_{j=1}^p
\dfrac{\P\big(\max (\sigma_{it},\sigma_{j0})\min (|Z_{it}|,|Z_{j0}|)>\gamma a_n\big)}{\P(|Z|>\gamma a_n)}\\
&\le &c\,\sum_{t=l}^{m_n} \sum_{i=1}^p\sum_{j=1}^p
\dfrac{\P\big(\sigma_{it}\min (|Z_{it}|,|Z_{j0}|)>\gamma a_n\big)}{\P(|Z|>\gamma a_n)}\,.
%+\dfrac{\P\big(\sigma_{0j}\min (|Z_{it}|,|Z_{0j}|)>\delta a_n\big)}{\P(|Z|>\delta a_n)}\Big]
\eeao
By stationarity the \pro ies on the \rhs\ do not depend on $t\ge l$. Therefore and by Breiman's result, the \rhs\ is bounded by
\beao
c\, m_n \dfrac{\P\big(\min (|Z_{it}|,|Z_{j0}|)>\gamma a_n\big)}{\P(|Z|>\gamma a_n)}=O((m_n/n) [n\,\,\P(|Z|>a_n)])=o(1)\,.
\eeao
This proves \eqref{eq:acl} for $(\bfX_t)$.\\[2mm]
(4)
Next we check the vanishing small values condition \eqref{eq:vansm} for the partial sums of 
$(\bfX_t^2)$ and $\alpha\in (2,4)$. 
It is not difficult to see that it suffices to prove the corresponding result for the component
processes:
\beam\label{eq:van2}
&&\lim_{\vep\downarrow 0}\limsup_{\nto}\P\Big(\Big|\sum_{t=1}^n \big(X_{it}^2\1(|X_{it}|\le \vep a_n)-\E [X_{it}^2\1(|X_{it}|\le \vep a_n)]\big)\Big|>\gamma a_n^2\Big)=0\,,\\&& \quad \gamma>0\,,\;i=1,\ldots,p\,.\nonumber
\eeam
%An application of a symmetrization inequality (Petrov \cite{petrov:1975}, Theorem 12 in Section III.3) conditional on $(\sigma_{it})$ 
%yields for an iid copy $(Z_{it}')$ of $(Z_{it})$: 
We have 
\beao\lefteqn{
a_n^{-2}\sum_{t=1} ^n\sigma_{it}^2
\E\big[Z_{it}^2\1(|X_{it}|\le \vep a_n)\mid \sigma_{it}] -a_n^{-2}\,n\,\E[X_{it}^2\1(|X_{it}|\le \vep a_n)]}\\
&=&  a_n^{-2}\sum_{t=1} ^n(\sigma_{it}^2 -\E[\sigma_{it}^2])\, \E[Z^2] -
a_n^{-2} \sum_{t=1} ^n \big(\sigma_{it}^2 \E[Z_{it}^2 \1(|X_{it}|>\vep a_n)\mid \sigma_{it}] - \E[X_{it}^2 \1(|X_{it}|>\vep a_n)]\big)\\&=&I_1+I_2\,.
\eeao
The \seq\ $(\sigma_{it}^2)$ satisfies the \clt\ with normalization $\sqrt{n}$. This follows from Ibragimov's \clt\ for strongly mixing \seq\ 
whose rate \fct\ $(\alpha_h)$ satisfies \eqref{eq:mixi} and has moment $\E [\sigma^{2(2+\delta))}]<\infty$ (see \eqref{eq:moment}); 
cf. Doukhan \cite{doukhan:1994}, p.~45. We know that $\sqrt{n}/a_n^2\to 0$ for $\alpha\in (2,4)$. Therefore $I_1\stp 0$.
We also have
\beao
\E[I_2^2]&\le& \,\dfrac n {a_n^{4}} \E \big[\sigma^4 (\E [Z^2\1(|X|>\vep a_n)\mid \sigma])^2\big]\\&& + 2\, \dfrac {n}{a_n^{4}} \sum_{h=1}^n |\cov(\sigma_{i0}^2 \E\big[Z_{i0}^2\1(|X_{i0}^2|> \vep a_n)\mid \sigma_{i0}], \sigma_{ih}^2 \E\big[Z_{ih}^2\1(|X_{ih}^2|> \vep a_n)\mid \sigma_{ih}])| \\
&=&I_3+I_4\,.
\eeao
%By Karamata's theorem 
In view of the moment conditions on $\sigma$ and since  $\E[Z^2]<\infty$,
$I_3\le c (n/a_n^4)\to 0$. In view of Doukhan \cite{doukhan:1994}, Theorem 3 on p.~9,  we have
\beao
I_4&\le& c\,\dfrac n {a_n^{4}}\sum_{h=1}^n \alpha_h^{\delta/(2+\delta)} (\E|\sigma|^{2(2+\delta)})^{2/(2+\delta)}\to 0\,.
\eeao
Thus it suffices for \eqref{eq:van2} to prove
\beao
\lim_{\vep\downarrow 0}\limsup_{\nto} \P\Big(\Big|\sum_{t=1}^n \big(\sigma_{it}^2\E[Z_{it}^2\1(|X_{it}|
\le \vep a_n)\mid\sigma_{it}] - X_{it}^2\1(|X_{it}|\le \vep a_n)\big)\Big|>\gamma\,a_n^2\Big)=0\,,\qquad \gamma>0\,.
\eeao
The summands are independent and centered, conditional on the $\sigma$-field generated by $(\sigma_{it})_{t=1,\ldots,n}$. An application of 
\v Cebyshev's inequality 
conditional on this $\sigma$-field and 
%independent centered random variable (see Petrov \cite{petrov:1995}, 2.6.20 on p.~82) 
Karamata's theorem yield, as $\nto$,
\beao&&
\E\Big[\P\Big(\Big|\sum_{t=1}^n \big(\sigma_{it}^2\E[Z_{it}^2\1(|X_{it}|
\le \vep a_n)\mid\sigma_{it}] - X_{it}^2\1(|X_{it}|\le \vep a_n)\big)\Big|>\gamma\,a_n^2\big| (\sigma_{is})\Big)\Big]\\
&\le &
c\, a_n^{-4}\E\Big[\sum_{t=1}^n \var(X_{it}^2\1(|X_{it}|\le \vep a_n)\mid \sigma_{it})\mid (\sigma_{is})\Big]\\
&\le &c\,n\,\vep^{4}\,
%\Big[\E\big[\sigma_{it}^{2\nu} \big[(\E[Z_{it}^2\1(|X_{it}|\le \vep a_n)\mid \sigma_{it}])^\nu+
\E[|X/(\vep a_n)|^{4}\1(|X|\le \vep a_n)]%\mid \sigma_{it}\big]\big]\Big]\\
\to  c\,\vep^{4-\alpha}\,.
\eeao
The \rhs\ converges to zero as $\vep\downarrow 0$.

This proves that all assumptions of Theorem \ref{thm:mikwin} are satisfied. 
Therefore the random variables on the \lhs\ of \eqref{eq:feller} converge to an $\alpha$-stable random vector with log-\chf
\begin{eqnarray*}
&&\int_0^\infty \E\big[\ex^{i\,y\,\bft'\sum_{j=0}^\infty \bfTh_j^{(2)}}- \ex^{i\,y\,\bft'\sum_{j=1}^\infty \bfTh_j^{(2)}}-i\,y\,\bft'\1_{(1,2)}(\alpha/2)\big]\,
d(-y^{\alpha/2})\\
&=&\sum_{j=1}^p \frac{1}{p}\int_0^\infty \E\big[\ex^{i\,y\,t_j}-i\,y\,t_j\1_{(1,2)}(\alpha/2)\big]\,
d(-y^{\alpha/2})\,,\qquad \bft=(t_1, \ldots, t_p)'\in\bbr^p,
\end{eqnarray*}
where we used \eqref{Eq:Theta_0:measure} and that $\bfTh_h^{(2)}=\bf0$ for $h\ge 1$. One easily checks that all summands in this expression are homogeneous functions in $t_j$ of degree $\alpha/2$.  Therefore, the limiting random vector in \eqref{eq:feller} has the same distribution as the sum $\sum_{j=1}^p \bfe_j \xi_{j,\alpha/2}$ for iid $\xi_{j,\alpha/2}$ which are $\alpha/2$-stable and totally skewed to the right (because all the summands in $S_j$ are non-negative).
\end{proof}

\subsection{Eigenvalues of the sample covariance matrix}\label{subsec:eigen}
We have the following approximations:
\bth\label{lem:1}
Assume that one of the following conditions holds:
\begin{enumerate}
\item
$(X_{it})$ is an iid field of \regvary\ \rv s with index 
$\alpha\in (0,4)$. If \ $\E[|X|]<\infty$ we also assume $\E[X]=0$.
\item
$(X_{it})$ is a \sv\ model \eqref{eq:1a} satisfying the \regvar , mixing and moment conditions of Theorem~\ref{the:1}.
If $\E [|Z|]<\infty$ we also assume $\E[Z]=0$. 
\end{enumerate}
Then, with $\mathbf{X}^n$ as in \eqref{eq:esses},
\beao
&&a_n^{-2} \twonorm{\bfX^n(\bfX^n)'-\diag(\bfX^n(\bfX^n)')}\stp 0\,,
\eeao
where $\twonorm{\cdot}$ is the spectral norm and 
$(a_n)$ is a \seq\ \st\ $n\,\P(|X|>a_n)\to 1$.
\ethe
\begin{proof} Part (1). Recall that for a 
$p\times p$ matrix $\bfA$ we have $\twonorm{\bfA}\le \frobnorm{\bfA}$, where $\frobnorm{\cdot}$ denotes the Frobenius norm.
Hence
\beam\label{eq:1}
a_n^{-4} \twonorm{\bfX^n(\bfX^n)'-\diag(\bfX^n(\bfX^n)')}^2&\le & a_n^{-4} \frobnorm{\bfX^n(\bfX^n)'-\diag(\bfX^n(\bfX^n)')}^2\nonumber \\
&=& \sum_{1\le i\ne j\le p} \big(a_n^{-2}S_{ij}\big)^2\,.
\eeam
In view of the assumptions, $(X_{it}\,X_{jt})_{t=1,2,\ldots}$, $i\ne j$, is an iid \seq\ of \regvary\ \rv s with index $\alpha$ which is 
also centered if $\E[|X|]<\infty$. We consider two different cases.\\
{\em The case $\alpha\in (0,2)$.} According to classical limit theory (see Feller \cite{feller:1971}, Petrov \cite{petrov:1995}) we have 
for $i\ne j$,
$
b_n^{-1}S_{ij}\std \xi_{\alpha}$, (see \eqref{eq:esses} for the definition of $S_{ij}$)
where $\xi_{\alpha}$ is an $\alpha$-stable \rv\ and $(b_n)$ is 
chosen \st\ $n\,\P(|X_1X_2|>b_n)\to 1$ for independent copies $X_1,X_2$ of $X$.
Since $(b_n)$ and $(a_n^2)$ are \regvary\ with indices 
$1/\alpha$ and $2/\alpha$, respectively, the \rhs\ in \eqref{eq:1} converges to zero in \pro y.\\[1mm]
{\em The case $\alpha\in [2,4)$.} In this case the \ds\ of $X_1X_2$ is in the domain of attraction of the normal law.
Since $X_1X_2$ has mean zero we can apply classical limit theory (see Feller \cite{feller:1971}, Petrov \cite{petrov:1995}) to conclude that
$
b_n^{-1} S_{ij} \std N\,,
$
where $(b_n)$ is \regvary\ with index $1/2$ and $N$ is centered Gaussian. Since $b_n/a_n^2\to 0$ we again conclude that the \rhs\ of \eqref{eq:1}
converges to zero in \pro y.\\[2mm]
Part (2). We again appeal to \eqref{eq:1}.
Let $\gamma < \min(2,\alpha)$. Then we have for $i \neq j$, using the independence of $(X_{it}X_{jt})$ conditional on $((\sigma_{it},\sigma_{jt}))$ and that the distribution of $Z$ is centered if its first absolute moments exists, that
\beao
a_n^{-2\gamma} \E\Big[\big|S_{ij}\big|^\gamma\mid ((\sigma_{it},\sigma_{jt}))\Big]
&\le &c\,\dfrac{n}{a_n^{2\gamma}} \dfrac 1 n \sum_{t=1}^n (\sigma_{it}\sigma_{jt})^{\gamma}(\E|Z|^\gamma)^2\, ,
\eeao
cf.\ von Bahr and Ess\'een \cite{bahr:esseen:1965} and Petrov~\cite{petrov:1995}, 2.6.20 on p.~82.
In view of the moment condition \eqref{eq:moment} we have 
$ \E[(\sigma_{i}\sigma_{j})^{\gamma}]<\infty$ and $n/a_n^{2\gamma}\to 0$ if we choose $\gamma$ sufficiently close to $\min(2,\alpha)$. Then the \rhs\ 
converges to zero in view of the ergodic theorem. An application of the conditional Markov inequality 
of order $\gamma$ yields 
$
a_n^{-2} S_{ij}\stp 0\,.
$
This proves the theorem.
\end{proof}
\bco\label{cor:sv1}
Assume that 
$(X_{it})$ is either
\begin{enumerate}
\item
an iid field of \regvary\ \rv s with index $\alpha\in (0,4)$ and $\E [X]=0$ if $\E [|X|]< \infty$, or
\item 
a \sv\ model of \regvary\ \rv s with index $\alpha\in (0,4)\setminus \{2\}$ satisfying the conditions of 
Theorem~\ref{lem:1}(2).
\end{enumerate}
Then
\beao
a_n^{-2}\max_{i=1,\ldots,p} \big|\la_{(i)}- S_{(i)}\big|\stp 0\,,
\eeao
where $(\la_i)$ are the eigenvalues of $\bfX^n(\bfX^n)'$, $\la_{(1)}\ge \cdots\ge \la_{(p)}$ are their ordered values  
and $S_{(1)}\ge \cdots \ge S_{(p)}$ are the ordered values of $S_1,\ldots,S_p$ defined in \eqref{eq:esses}. 
In particular, we have 
\beam\label{eq:2a}
a_n^{-2} \big(\la_{(1)}-c_n,\ldots,\la_{(p)}-c_n\big) \std \big(\xi_{(1),\alpha/2},\ldots,\xi_{(p),\alpha/2}\big)\,,
\eeam
where $(c_n)$ is defined in \eqref{eq:6} for $\alpha\ne 2$ and in Remark~\ref{rem:feller} for $\alpha=2$,
$(\xi_{i,\alpha/2})$ are iid $\alpha/2$-stable \rv s given in Theorem~\ref{the:1} for the \sv\ model and in Remark~\ref{rem:feller} for the iid field, 
and $\xi_{(1),\alpha/2}\ge \cdots\ge \xi_{(p),\alpha/2}$ are their
ordered values. 
\eco
\begin{proof}
We have by Weyl's inequality (see Bhatia \cite{bhatia:1997}) and  Theorem~\ref{lem:1},
\beam\label{eq:weyl1}
a_n^{-2}\max_{i=1,\ldots,p} \big|\la_{(i)}- S_{(i)}\big|&\le &
a_n^{-2} \twonorm{\bfX^n(\bfX^n)'-\diag(\bfX^n(\bfX^n)')}\stp 0\,.
\eeam
If $(X_{it})$ is an iid random field (see Remark~\ref{rem:feller}) or a \sv\ model satisfying the conditions of Theorem~\ref{lem:1}(2)
we have \eqref{eq:feller}. Then \eqref{eq:weyl1} implies \eqref{eq:2a}.
\end{proof}
\bre\label{rm:many}\rm 
If $\alpha\in (2,4)$ we have $\E[X^2]<\infty$. Therefore
\eqref{eq:2a} reads as 
\beam\label{eq:rem:many}
\dfrac{n}{a_n^2}\,\big(\dfrac{\la_{(i)}}{n}- \E[X^2]\big)_{i=1,\ldots,p}\std (\xi_{(i),\alpha/2})_{i=1,\ldots,p}\,.
\eeam
We notice that $n/a_n^2\to \infty$ for $\alpha\in (2,4)$ since $(n/a_n^2)$ is \regvary\ with index $1-2/\alpha$. In particular, if ${\rm tr}(\bfX^n(\bfX^n)')$ denotes the trace of $\bfX^n(\bfX^n)'$ we have for $i\le p$,
\beam\label{eq:ratio}
\dfrac{\la_{(i)}}{{\rm tr}(\bfX^n(\bfX^n)')}&=& \dfrac{\la_{(i)}/n}{(\la_1+\cdots +\la_p)/n}
\stp \dfrac  1p\,.
%\dfrac{\la_{(1)}+\cdots + \la_{(i)}}{{\rm tr}(\bfX\bfX')}&\stp& \dfrac  ip\,.
\eeam
\par
The joint \asy\ distribution of the ordered eigenvalues $(\la_{(i)})$ is easily calculated from the distribution of a totally skewed
$\alpha/2$-stable \rv\ $\xi_{1,\alpha/2}$; in particular, the limit of $(a_n^{-2}(\la_{(1)}-c_n))$ has the \ds\ of $\max(\xi_{1,\alpha/2},\ldots,\xi_{p,\alpha/2})$.
\par
For applications, it is more natural to replace the \rv s $X_{it}$ by their mean-centered versions $X_{it}-\ov X_i$,
where $\ov X_{i}= (1/n) \sum_{t=1}^n X_{it}$, instead of assuming that they have mean zero. The previous results remain valid
for the sample-mean centered \rv s $X_{it}$, also in the case when $X$ has infinite first moment.
\ere

\subsection{Some applications: Limit results for ordered eigenvalues and eigenvectors of the sample covariance matrix}\label{subsec:appl}
In what follows, we assume the conditions of Corollary~\ref{cor:sv1}.
\subsubsection{Spacings}\label{ssection:spacing}
Using the joint \con\ of the normalized ordered eigenvalues $(\la_{(i)})$ we can calculate the limit of the spectral gaps:
\beam\label{eq:order}
\big(\dfrac{\la_{(i)}-\la_{(i+1)}}{a_n^2}\big)_{i=1,\ldots,p-1} \std \big(\xi_{(i),\alpha/2}-\xi_{(i+1),\alpha/2}\big)_{i=1,\ldots,p-1}\,.
\eeam
\par
We notice that the ordered values $\xi_{(i),\alpha/2}$ and linear \fct als thereof (such as $\xi_{(i),\alpha/2}-\xi_{(i+1),\alpha/2}$) are again jointly \regvary\ with index $\alpha/2$.
This is due to the continuous mapping theorem for \regvary\ vectors; see Hult and 
Lindskog \cite{hult:lindskog:2005,hult:lindskog:2006}, cf. 
Jessen and Mikosch \cite{jessen:mikosch:2006}.
\subsubsection{Trace}
For the trace of $\bfX^n(\bfX^n)'$ we have
\beao
a_n^{-2} \big({\rm tr}(\bfX^n(\bfX^n)')-p\,c_n\big)&=& a_{n}^{-2}\sum_{i=1}^p (S_i-c_n)\\
&=& a_{n}^{-2}\sum_{i=1}^p (\la_i-c_n)\std \xi_{1,\alpha/2}+\cdots+\xi_{p,\alpha/2}\eqd p^{2/\alpha}\xi_{1,\alpha/2}\,.
\eeao
Moreover, we have the joint \con\ of the normalized and centered $(\la_{(i)})$ and ${\rm tr}(\bfX^n(\bfX^n)')=\la_1+\cdots +\la_p$. In particular,
we have the self-normalized limit relations
\beao
\big(\dfrac{\la_{(i)}-c_n}{{\rm tr}(\bfX^n(\bfX^n)')-p\,c_n}\big)_{i=1,\ldots,p}&\std& 
\big(\dfrac{\xi_{(i),\alpha/2}}{\xi_{1,\alpha/2}+\cdots+\xi_{p,\alpha/2}}\big)_{i=1,\ldots,p}\,,
\eeao
and for $\alpha\in (2,4)$, by the \slln ,
\beao
\dfrac{np}{a_n^2}\,\Big(\dfrac{\la_{(i)}-c_n}{{\rm tr}(\bfX^n(\bfX^n)')}\Big)_{i=1,\ldots,p}&\std& \dfrac{\xi_{(i),\alpha/2}}{\E [X^2]}\,.
\eeao
\subsubsection{Determinant}
Since $\la_i-c_n$ are the eigenvalues of $\bfX^n(\bfX^n)'-c_n \bfI_p$, where $\bfI_p$ is the $p\times p$ identity matrix,
we obtain for the determinant
\beao
{\rm det} \big(a_n^{-2}(\bfX^n(\bfX^n)'-c_n \,\bfI_p)\big)&=& \prod_{i=1}^p a_n^{-2} (\la_{(i)}-c_n)\\
&\std& \xi_{(1),\alpha/2}\cdots \xi_{(p),\alpha/2}=\xi_{1,\alpha/2}\cdots \xi_{p,\alpha/2} \,. 
\eeao
For $\alpha\in (2,4)$, we also have 
\beao
\dfrac 1{a_n^2 c_n^{p-1}} \big({\rm det} (\bfX^n(\bfX^n)')- c_n^p\big)&=&
\sum_{i=1}^p  a_n^{-2}\big(\la_{(i)}-c_n\big) \prod_{j=1}^{i-1} \dfrac{\la_{(j)}}{c_n}\\
&\std &\sum_{i=1}^p \xi_{(i),\alpha/2}= \sum_{i=1}^p\xi_{i,\alpha/2} \eqd p^{2/\alpha}\,\xi_{1,\alpha/2}\,,
\eeao
where we used \eqref{eq:rem:many}. 
\subsubsection{Eigenvectors}\label{Subsub:eigenvectors}
It is also possible to localize the eigenvectors of the matrix $a_n^{-2}\bfX^n(\bfX^n)'$. Since this matrix is approximated
by its diagonal in spectral norm, one may expect that the unit eigenvectors of the original matrix are close to the canonical basis
vectors. We can write
\beao
a_n^{-2}\bfX^n(\bfX^n)'\bfe_{L_j}= a_n^{-2}S_{(j)}\,\bfe_{L_j}+ \vep_n\,\bfW\,,
\eeao
where $\bfW$ is a unit vector orthogonal to $\bfe_{L_j}$, $L_j$ is the index of $S_{(j)}=S_{L_j}$ and 
\beao
\vep_n= a_n^{-2}\|\big(\bfX^n(\bfX^n)'- S_{(j)}\big) \bfe_{L_j}\|_{\ell_2}\stp 0\,,
\eeao
from Theorem \ref{lem:1} and by equivalence of all matrix norms. 
According to Proposition A.1 in Benaych-Georges and Pech\'e \cite{benaych:peche:2014}, there is an eigenvalue $a_n^{-2}\la_{(j)}$  
of  $a_n^{-2}\bfX^n(\bfX^n)'$ in some $\vep_n$-neighborhood of  $a_n^{-2}S_{(j)}$. Define
\beao\Omega_n=\{a_n^{-2}|\la_{(j)}-\la_{(l)}|>d_n\,,l\ne j \}\,,
\eeao
for $d_n=k\vep_n$ for any fixed $k>1$. Then
$\lim_{\nto}\P(\Omega_n)=1$ because of \eqref{eq:order} and $d_n\stp 0$.
Hence, for large $n$, $a_n^{-2}\la_{(j)}$ and $a_n^{-2}\la_{(l)}$ have distance at least $d_n$ with high \pro y. Another application of
Proposition A.1 in \cite{benaych:peche:2014} yields that the unit eigenvector $\bfV$ associated with $a_n^{-2}\la_{(j)}$ 
satisfies the relation
\beao
\limsup_{\nto}\P\big(\|\bfV- V_{L_j}\bfe_{L_j}\|_{\ell_2}>\delta\big)&\le &
\limsup_{\nto}\P\big(\{\|\bfV- V_{L_j}\bfe_{L_j}\|_{\ell_2}>\delta\}\cap \Omega_n\big)+
\limsup_{\nto}\P(\Omega_n^c)\\
&\le &\limsup_{\nto}\P\big( \{2\,\vep_n /(d_n-\vep_n)>\delta\}\cap \Omega_n\big)\\
%& \le & \limsup_{\nto}\P\big( 2\,\vep /(d-\vep)>\delta\big)\\
&= & \1_{\{2/(k-1)>\delta\}}\,.
\eeao
For any fixed $\delta>0$, the \rhs\ is zero for sufficiently large $k$. Since both $\bfV$ and $\bfe_{L_j}$ are unit eigenvectors
this means that 
$\|\bfV- \bfe_{L_j}\|_{\ell_2}\stp 0$.
\subsubsection{Sample correlation matrix}
In Remark~\ref{rm:many} we mentioned that we can replace the variables $X_{it}$ by their sample-mean centered versions $X_{it}-\ov X_i$
without changing the \asy\ theory. Similarly, one may be interested in transforming the $X_{it}$ as follows:
\beao
\wt X_{it}= \dfrac{X_{it}- \ov X_i}{\wh \sigma_i}\,,\qquad 
\wh \sigma_i^2= \sum_{t=1}^n (X_{it}-\ov X_i)^2\,. 
\eeao
Then the matrix
\beao
\wt\bfX^n(\wt \bfX^n)'= \big(\sum_{t=1}^n\wt X_{it}\wt X_{jt}\big)_{i,j=1,\ldots,p}\,,
\eeao
is  the sample correlation matrix. We write $\wt \la_i$, $i=1,\ldots,p$, for the eigenvalues
of $\wt \bfX^n(\wt \bfX^n)'$ and $\wt \la_{(1)}\ge \cdots \ge \wt \la_{(p)}$ for their ordered values. 
\par
We notice that the entries of this matrix are all bounded in modulus by one.
In particular, the diagonal consists of ones.
We do not have a complete  limit theory for the eigenvalues  $\wt \la_i$.
We restrict ourselves to iid $(X_{it})$ to explain the differences.
\begin{lemma}
Assume that 
$(X_{it})$ is an iid field of \rv s. 
%If $\E[|X|]<\infty$ we also assume $\E[X]=0$.
\begin{enumerate}
\item If $\E[X^2]<\infty$ then 
\beao
\sqrt{n} \max_{i=1,\ldots,p} |\wt \la_{i}-1|= O_\P(1)\,.
\eeao
\item If $X$ is \regvary\ with index $\alpha\in (0,2)$ then
\beao
\dfrac{a_n^2}{b_n} \max_{i=1,\ldots,p} |\wt \la_{i}-1|=O_\P(1)\,,
\eeao
where $(a_n)$ and $(b_n)$  are chosen \st\ $\P(|X|>a_n)\sim \P(|X_1X_2|>b_n) \sim n^{-1}$ for iid copies $X_1,X_2$ of $X$. 
\end{enumerate}
\end{lemma}
\bre
Notice that the lemma implies $\wt \la_i\stp 1$ for $i=1,\ldots,p$, and the analog of
relation \eqref{eq:ratio} remains valid. 
\ere
\begin{proof} Part(1) We assume without loss of generality that $1=\E [X^2]$.
Then by classical limit theory,
\beao
\sqrt{n}\big(\wt \bfX^n(\wt \bfX^n)'-{\rm diag} (\wt \bfX^n(\wt \bfX^n)')\big)&=& \sqrt{n}\big(\wt \bfX^n(\wt \bfX^n)'-\bfI_p\big)\\
&=& \Big(\1(i\ne j)\dfrac{n^{-1/2}\sum_{t=1}^n (X_{it}-\ov X_{i})(X_{jt}-\ov X_j)}{(\wh \sigma_i/\sqrt{n}) 
(\wh\sigma_j/\sqrt{n})}\Big)\\
&\std & \big(N_{ij}\1(i\ne j)\big)\,,
\eeao
where $N_{ij}$, $1\le i<j\le n$, are iid $N(0,1)$ and $N_{ij}=N_{ji}$. By Weyl's inequality,
\beao
\sqrt{n} \max_{i=1,\ldots,p}\Big|\wt \la_{(i)}-1\Big|\le\sqrt{n} \|\wt \bfX^n(\wt \bfX^n)'-\bfI_p\|_2=O_\P(1)\,.
\eeao
Part(2) If $X$ is \regvary\ with index $\alpha\in (0,2)$, we have that $(a_n^{-2} \wh \sigma_i^2)$ converges to a vector of iid positive $\alpha/2$-stable
\rv s $(\xi_i)$, while for every $i\ne j $, $b_n^{-1}\sum_{t=1}^n 
(X_{it}-\ov X_{i})\,(X_{jt}-\ov X_{j})\std \xi_{ij}$ and the limit $\xi_{ij}$ is $\alpha$-stable.
Then by Weyl's inequality
\beao
\dfrac{a_n^2}{b_n} \max_{i=1,\ldots,p}\Big|\wt \la_{(i)}-1\Big|\le \dfrac{a_n^2}{b_n}\|\wt \bfX^n(\wt \bfX^n)'-\bfI_p\|_2=O_\P(1)\,.
\eeao
\end{proof}
\section{Case (2): $\sigma$ dominates the tail}\label{Sec:case2}\setcounter{equation}{0}
In this section we assume the conditions of Case (2); see Section \ref{subsec:Case:2}. Our goal is to derive results analogous to Case (1): \regvar\ 
of $(X_{it})$, infinite variance limits for $S_{ij}$ and limit theory for the eigenvalues of the corresponding sample 
covariance matrices. It turns out that this case offers a wider spectrum of possible limit behaviors and that we have to further distinguish our assumptions about the distribution of $\eta$. So, in addition to \eqref{Eq:tail:eta} we assume that either
\begin{equation}\label{Eq:inf:mom:eta} \E[\ex^{\eta \alpha}]=\infty\end{equation}
or 
\begin{equation}\label{Eq:conv:eq:eta} \lim_{x \to \infty}\frac{\P(\eta_1+\eta_2>x)}{\P(\eta_1>x)}=c \in (0,\infty) \;\; \Leftrightarrow\;\; \lim_{x \to \infty}\frac{\P(\ex^{\eta_1}\cdot \ex^{\eta_2}>x)}{\P(\ex^{\eta_1}>x)}=c \in (0,\infty)\end{equation}
hold, where $\eta_1$ and $\eta_2$ are independent copies of $\eta$.
\bre Following  Cline \cite{cline:1986}, we call 
the distribution of a \rv\ $\eta$ {\it convolution equivalent} if  $\ex^\eta$ is regularly varying and
relation \eqref{Eq:conv:eq:eta} holds. The assumptions \eqref{Eq:inf:mom:eta} and \eqref{Eq:conv:eq:eta} are mutually exclusive, since the only possible finite limit $c$ in \eqref{Eq:conv:eq:eta} is given by $c=2\E[\ex^{\eta \alpha}]$; see Davis and Resnick \cite{davis:resnick:1986}. There are, however, 
regularly varying distributions of $\ex^\eta$ which satisfy $\E[\ex^{\eta \alpha}]<\infty$ but not
\eqref{Eq:conv:eq:eta}. An example is given in Cline \cite{cline:1986}, p.\ 538;  see also Lemma~\ref{lem:product}(3) 
for a necessary and sufficient condition ensuring \eqref{Eq:conv:eq:eta}.
\ere
As we will see later, relations \eqref{Eq:inf:mom:eta} and \eqref{Eq:conv:eq:eta} cause
rather distinct limit behavior of the sample covariance matrix. In particular, \eqref{Eq:conv:eq:eta}
allows for non-vanishing off-diagonal elements of the normalized sample  covariance matrices, in contrast to Case (1). 
\par
For notational simplicity, define
\beao
\psi=\max_{k,l} \psi_{kl}\qquad \mbox{and}\qquad \Lambda=\{(k,l): \psi_{kl}=\psi\}\,.
\eeao
Recall that for convenience we assume that $\psi=1$; if the latter condition does not hold we can replace  
(without loss of generality) the \rv s $\eta_{kl}$ by $\psi \eta_{kl}$ and the coefficients $\psi_{kl}$ by $\psi_{kl}/\psi$.
For given $(i,j)$, we define 
\begin{equation} \label{eq:psi:ij:def}
\psi^{ij}= \max_{k,l}\, (\psi_{kl}+\psi_{k+i-j,l})\,.
\end{equation}
Notice that $1\le \psi^{ij}\le 2$. For $d\ge 1$, we write $\mathbf{i}=(i_1, \ldots, i_d), \mathbf{j}=(j_1, \ldots, j_d)$ for elements of $\mathbb{Z}^d$. For given  $\mathbf{i}$ and  $\mathbf{j}$ we also define
\beao
\psi^{\mathbf{i},\mathbf{j}}= \max_{1 \leq l \leq d} \psi^{i_l,j_l}.
\eeao

\subsection{Regular variation}\label{subsec:regvar2}
We start by showing that the volatility sequences are regularly varying.
\bpr\label{prop:rv:1} Under the aforementioned conditions and conventions (including that 
either \eqref{Eq:inf:mom:eta} or \eqref{Eq:conv:eq:eta} hold), 
\begin{enumerate}
\item
each of the \seq s 
$(\sigma_{it})_{t \in \mathbb{Z}}$, $i=1,2,\ldots$, is \regvary\ with index $\alpha$,
\item
each of the \seq s $(\sigma_{it}\sigma_{jt})_{t \in \mathbb{Z}}$, $i,j=1,2,\ldots$, is \regvary\ with corresponding index  $\alpha/\psi^{ij}$,
\item
For $d\geq 1$ and $\mathbf{i},\mathbf{j}\in\bbz$, the $d$-variate \seq\ 
$((\sigma_{i_k,t}\sigma_{j_k,t})_{1\leq k \leq d})_{t \in \mathbb{Z}}$ is regularly 
varying with index $\alpha/\psi^{\mathbf{i},\mathbf{j}}$.
\end{enumerate}
\epr
\bre Part (3) of the proposition possibly includes degenerate cases in the sense that for some choices of $(i_k,j_k)$, 
$(\sigma_{i_k,t}\sigma_{j_k,t})$ is \regvary\ with index $\alpha/\psi^{i_k,j_k}> \alpha/\psi^{\mathbf{i},\mathbf{j}}$. 
\par
Part (3) implies (2) in the case $d=1$. Part (2) implies (1) by setting $i=j$ and observing that, 
by non-negativity of $\sigma$,  \regvar\ of $(\sigma_{it}^2)$
with index $\alpha/2$ is equivalent to \regvar\ of $(\sigma_{it})$ with index $\alpha$. 
\ere
% \bre We notice that $(\sigma_{it}^2)_{t=1,2,\ldots}$ is \regvary\ with index $\alpha/(2\psi)$. If $m(i,j)>2$ for some $i\ne j$ then
% $(\sigma_{it}\sigma_{jt})_{t=1,2,\ldots}$ is \regvary\ with a smaller index than $(\sigma_{it}^2)_{t=1,2,\ldots}$. If $\Lambda$ consists of exactly one element
% then $m(i,j)=1$ for $i\ne j$.
% \ere
\begin{proof} To give some intuition we start with the proof of the marginal regular variation of $\sigma$, although it is just a special case of (1).
We have
\beam\label{eq:oo}
\sigma_{it}= \ex^{ \sum_{(k,l)\in\Lambda}\eta_{i-k,t-l}}\,
\ex^{\sum _{(k,l)\not \in\Lambda} \psi_{kl} \eta_{i-k,t-l}}=:
\sigma_{it,\Lambda}\sigma_{it,\Lambda^c}.
\eeam
%{\blue where we note that $\sigma_{it,\Lambda}$ and $\sigma_{it,\Lambda^c}$ are independent for each $i,t$.}
We first verify that $\sigma=\sigma_\Lambda\sigma_{\Lambda^c}$ is \regvary\ with index 
$\alpha$. Since $|\Lambda|<\infty$ by our assumptions, and in view of Embrechts and Goldie \cite{embrechts:goldie:1980}, Corollary on p.~245, cf. also Lemma~\ref{lem:product}(1) below,
the product
$\sigma_\Lambda$ is \regvary\ with index $\alpha$. % hence $\sigma_\Lambda$ is
%\regvary\ with index $\alpha/\psi$.
The \rv\ $\sigma_{\Lambda^c}$ is independent of $\sigma_\Lambda$. Similarly to
Mikosch and Rezapour \cite{mikosch:rezapour:2013} (see also the end of this proof for a similar argumentation) one can show that 
$\sigma_{\Lambda^c}$ has moment of order $\alpha+\vep$ 
for sufficiently small positive $\vep$. Therefore, by Breiman's lemma \cite{breiman:1965},
\beao
\P(\sigma >x)\sim \E [\sigma_{\Lambda^c}^{\alpha}]\,\P(\sigma_\Lambda>x)\,,\qquad \xto\,.
\eeao
This proves \regvar\ with index $\alpha$ of the marginal \ds s of $(\sigma_{it})$.
\par
In the remainder of the proof we focus on (3).
%Statement (2) is then a special case of (3) with $d=1$ and statement (1) follows by the non-negativity of $\sigma$ from (2) by setting $i=j$. 
%In the following let $d\geq 1$ and set $\mathbf{i}=\{i_1, \ldots, i_d\}, \mathbf{j}=\{j_1, \ldots, j_d\},$ $\mathbf{t}=\{t_1, \ldots, t_d\}$ all be subsets of $\mathbb{Z}$. Define
% $$ \psi^{\mathbf{i},\mathbf{j}}= \max_{1 \leq l \leq d} \psi^{i_l,j_l}.$$
For a given choice of $\mathbf{i, j, t}\in \bbz^d$, we write
% $$ \Lambda_{\mathbf{i},\mathbf{j},\mathbf{t}}^{(l)}=\{m,n: \psi_{m+i_l,n+t_l}+\psi_{m+j_l,n+t_l}= \psi^{\mathbf{i},\mathbf{j}}\}, \;\;\; 1 \leq l \leq d, $$
%and 
\begin{equation}\label{Eq:Lambda:set:def} \Lambda_{\mathbf{i},\mathbf{j},\mathbf{t}}=\{(m,n ): \psi_{i_l-m,t_l-n}+\psi_{j_l-m,t_l-n}= \psi^{\mathbf{i},\mathbf{j}}\; \mbox{for some } 1\leq l \leq d\}.\end{equation}
We will show that the random vector $(\sigma_{i_1,t_1}\sigma_{j_1,t_1}, \ldots, \sigma_{i_d,t_d}\sigma_{j_d,t_d})=:\boldsymbol{\sigma}'$ is regularly varying with index $\alpha/\psi^{\mathbf{i},\mathbf{j}}$ which proves (3). Note that
\begin{eqnarray*}
\sigma_{i,t}\sigma_{j,t}&=&\prod_{(k,l)}\exp(\psi_{kl}\eta_{i-k,t-l})\prod_{(k',l' )}\exp(\psi_{k'l'}\eta_{j-k',t-l'}) \\
&=& \prod_{(m,n )}\exp((\psi_{i-m,t-n}+\psi_{j-m,t-n})\eta_{m,n})
\end{eqnarray*}
and write
\begin{align} \nonumber \boldsymbol{\sigma}&= \mbox{diag}\left(\left(
  \begin{matrix}
  \prod\limits_{(m,n) \in \Lambda_{\mathbf{i},\mathbf{j},\mathbf{t}}^c}\ex^{\eta_{m,n}(\psi_{i_1-m,t_1-n}+\psi_{j_1-m,t_1-n})} \\ \vdots \\ \prod\limits_{(m,n) \in \Lambda_{\mathbf{i},\mathbf{j},\mathbf{t}}^c}\ex^{\eta_{m,n}(\psi_{i_d-m,t_d-n}+\psi_{j_d-m,t_d-n})}
  \end{matrix}
\right)'\right)\cdot \left(\begin{matrix}
                   \prod\limits_{(m,n) \in \Lambda_{\mathbf{i},\mathbf{j},\mathbf{t}}}\ex^{\eta_{m,n}(\psi_{i_1-m,t_1-n}+\psi_{j_1-m,t_1-n})} \\ \vdots \\ \prod\limits_{(m,n) \in \Lambda_{\mathbf{i},\mathbf{j},\mathbf{t}}}\ex^{\eta_{m,n}(\psi_{i_d-m,t_d-n}+\psi_{j_d-m,t_d-n})}
                   \end{matrix}
 \right)\\
\label{Eq:factors:sigma} &=: \mathbf{A}\, \mathbf{Z},
\end{align}
where ${\rm diag}((a_1,\ldots,a_k))$ is any diagonal matrix with diagonal elements $a_1,\ldots,a_k$. 
We notice that $\bfA$ and $\bfZ$ are independent.
\par
Consider iid copies $(Y_j)$ of $\ex^\eta$. There exist suitable numbers $(a_{ij})_{1\le i\le d,1\le j\le p}$ with $p= |\Lambda_{\mathbf{i},\mathbf{j},\mathbf{t}}|$
\st\ the components of $\mathbf{Z}$ have \rep\ in \ds\ 
$\prod_{j=1}^{p}Y_j^{a_{ij}}$, $1 \leq i \leq d $. By assumption, $Y_j$ is regularly varying with index $\alpha$ and
satisfies either assumption \eqref{Eq:conv:equivalent} or $\E[Y_j^\alpha]=\infty$. 
Furthermore, for each $j$ there exists one $1\leq i \leq d$ 
such that $a_{ij}=a_{\max}=\psi^{\mathbf{i},\mathbf{j}}$ by the definition of $\Lambda_{\mathbf{i},\mathbf{j},\mathbf{t}}$. 
An application of Proposition \ref{Pr:genRVforproducts} shows that $\bfZ$ is \regvary\ 
with index $\alpha/\psi^{\mathbf{i},\mathbf{j}}$ and limit measure $\mu_{\mathbf{Z}}$ which is given as $\mu$ in Proposition \ref{Pr:genRVforproducts} (ii) (if \eqref{Eq:inf:mom:eta} holds) or Proposition \ref{Pr:genRVforproducts} (i) (if \eqref{Eq:conv:eq:eta} holds). Now, choose $\epsilon, \delta>0$ such that
$$\frac{\psi_{i_l-m,t_l-n}+\psi_{j_l-m,t_l-n}}{\psi^{\mathbf{i},\mathbf{j}}}(1+\delta)<1-\epsilon, \;\;\;  (m,n) \in \Lambda_{\mathbf{i},\mathbf{j},\mathbf{t}}^c,\;\;1\leq l \leq d,$$
which is possible by the definition of $\Lambda_{\mathbf{i},\mathbf{j},\mathbf{t}}$ and the summability constraint on the coefficients. Then we have 
\begin{eqnarray*}&&\E\left[\|\mathbf{A}\|_{\mbox{\scriptsize op}}^{\alpha(1+\delta)/\psi^{\mathbf{i},\mathbf{j}}}\right]\\
&\leq & \sum_{l=1}^d \prod\limits_{(m,n) \in \Lambda_{\mathbf{i},\mathbf{j},\mathbf{t}}^c} \E\left[\ex^{\eta_{m,n}\alpha(1+\delta)(\psi_{i_l-m,t_l-n}+\psi_{j_l-m,t_l-n})/\psi^{\mathbf{i},\mathbf{j}}}\right]\\
&\leq& \sum_{l=1}^d \prod\limits_{(m,n) \in \Lambda_{\mathbf{i},\mathbf{j},\mathbf{t}}^c} \E\left[\ex^{\eta_{m,n}\alpha(1-\epsilon)}\right]^{(1+\delta)(\psi_{i_l-m,t_l-n}+\psi_{j_l-m,t_l-n})/((1-\epsilon)\psi^{\mathbf{i},\mathbf{j}})} <\infty,
\end{eqnarray*}
where we used Jensen's inequality for the penultimate step and the 
summability condition of the coefficients for the final one. Thus we have verified all conditions of the 
multivariate Breiman lemma in Basrak et al.~\cite{basrak:davis:mikosch:2002a}, implying that
$\boldsymbol{\sigma}$ inherits \regvar\  from $\bfZ$ with corresponding
index $\alpha/\psi^{\mathbf{i},\mathbf{j}}$ and limit measure
$\mu_{\boldsymbol{\sigma}}(\cdot)=\E[\mu_{\mathbf{Z}}(\mathbf{A}^{-1}\cdot)]$.
%with $\mathbf{A}$ as in \eqref{Eq:factors:sigma}.
\end{proof}
\bpr\label{prop:rv:2} Assume that the aforementioned conditions (including either \eqref{Eq:inf:mom:eta} or \eqref{Eq:conv:eq:eta}) hold and that in addition $\E[|Z|^{\alpha+\delta}]<\infty$ for some $\delta>0$. Then the following statements hold: 
\begin{enumerate}
\item
Each of the \seq s 
$(X_{it})_{t \in \mathbb{Z}}$, $i \in \mathbb{Z}$, is \regvary\ with index $\alpha$.\\[1mm]
If  \eqref{Eq:inf:mom:eta} holds then the corresponding spectral 
tail process satisfies $\Theta_t^{i}=0$ a.s., $t\ge 1$, and
$ \P(\Theta^{i}_0=\pm 1)=\E[Z_{\pm}^{\alpha}]/\E[|Z|^{\alpha}]$.\\[1mm]
 %\; \P(\theta_0^{(i)}=-1)=\frac{\E[(Z^-)^{\alpha}]}{\E[|Z|^{\alpha}]}
%\qquad \bfTh_h^{i}=0 \;\as,  h \geq 1\,.
If \eqref{Eq:conv:eq:eta} holds, then for any Borel set $B=B_0\times\cdots \times B_n\subset \bbr^{n+1}$,
\begin{align}\label{Eq:theta:law:1} \P((\Theta^{i}_t)_{t=0,\ldots,n} \in B)=\sum_{(u,v) \in \Lambda_i^{(0)}}\frac{1}{|\Lambda_i^{(0)} |}\frac{\E\left[\1\left(\left(\1((u,v) \in \Lambda_i^{(t)})\frac{X_{it}}{|X_{i0}|}\right)_{t=0, \ldots, n} \in B \right)|X_{i0}|^{\alpha}\right]}{\E[|X_{i0}|^{\alpha}]},
\end{align}
where $\Lambda_i^{(t)}=\{(u,v):\psi_{i-u,t-v}=1 \}, t=0, \ldots, n$. 
%the random variables $\ex^{\eta_{u,v}^{(u,v)}}$,$(\ex^{\eta_{l,m}^{(u,v)}})_{(l,m)\neq(u,v)}$, $(\tilde{Z}_{i,t})_{t=0, \ldots, n}$ 
%have joint {\red don't understand} $\nu_{\alpha} \otimes (\otimes_{(l,m)\neq(u,v)} P^{\ex^{\eta}}) \otimes (\otimes_{t=0, \ldots, n}P^{Z})$-density 
%$$ f^{(u,v)}(y_{u,v},(y_{l,m})_{(l,m) \neq (u,v)}, (z_{i,t})_{t=0, \ldots, n})=\frac{\1(\prod_{(l,m)}y_{l,m}^{\psi_{i-l,-m}}|z_{i,0}|>1)}{\E[(\prod_{(l,m)\neq(u,v)}\ex^{\eta_{l,m}^{(u,v)}\psi_{i-l,-m}}|\tilde{Z}_{i,0}|)^{\alpha}]},$$
\item
Each of the \seq s $(X_{it}X_{jt})_{t \in \mathbb{Z}}$, $i,j\in \mathbb{Z}$, is \regvary\ with index  $\alpha/\psi^{ij}$.\\[1mm] 
If \eqref{Eq:inf:mom:eta} holds then the corresponding spectral tail process satisfies $\Theta_t^{ij}=0$ a.s., $t\ge 1$, and
$\P(\Theta^{ij}_0=\pm 1)=\E[(Z_iZ_j)_{\pm}^{\alpha/\psi^{ij}}]/\E[|Z_iZ_j|^{\alpha/\psi^{ij}}]$.\\[1mm]
%, \; \P(\theta_0^{(ij)}=-1)=\frac{\E[((Z_iZ_j)^-)^{\alpha/\psi^{ij}}]}{\E[|Z_iZ_j|^{\alpha/\psi^{ij}}]}, $$
If \eqref{Eq:conv:eq:eta} holds, then for any Borel set $B=B_0\times\cdots \times B_n\subset \bbr^{n+1}$, 
\begin{align}  
&\P(( \Theta^{ij}_t)_{t=0,\ldots,n} \in B) \nonumber \\
 \label{Eq:theta:law:2}=&\sum_{(u,v) \in \Lambda_{i,j}^{(0)}}\frac{1}{|\Lambda_{i,j}^{(0)} |}\frac{\E\left[\1\left(\left(\1((u,v) \in \Lambda_{i,j}^{(t)})\frac{X_{it}X_{jt}}{|X_{i0}X_{j0}|}\right)_{t=0, \ldots, n} \in B \right)|X_{i0}X_{j0}|^{\alpha/\psi^{ij}}\right]}{\E[|X_{i0}X_{j0}|^{\alpha/\psi^{ij}}]},  
\end{align}
where $\Lambda_{i,j}^{(t)}=\{(u,v):\psi_{i-u,t-v}+\psi_{j-u,t-v}=\psi^{ij}\}, t=0, \ldots, n$. 
%and the random variables $\ex^{\eta_{u,v}^{(u,v)}}$,$(\ex^{\eta_{l,m}^{(u,v)}})_{(l,m)\neq(u,v)}$, $(\tilde{Z}_{i,t},\tilde{Z}_{j,t})_{t=0, \ldots, n}$ have joint $\nu_{\alpha} \otimes (\otimes_{(l,m)\neq(u,v)} P^{\ex^{\eta}}) \otimes (\otimes_{t=0, \ldots, n}P^{(Z_{i},Z_{j})})$-density 
%\begin{eqnarray*}&& f^{(u,v)}(y_{u,v},(y_{l,m})_{(l,m) \neq (u,v)}, (z_{i,t},z_{j,t})_{t=0, \ldots, n})\\
%&=&\frac{\1(\prod_{(l,m)}y_{l,m}^{\psi_{i-l,-m}+\psi_{j-l,-m}}|z_{i,0}z_{j,0}|>1)}{\E[(\prod_{(l,m)\neq(u,v)}\ex^{\eta_{l,m}^{(u,v)}(\psi_{i-l,-m}+\psi_{j-l,-m})}|\tilde{Z}_{i,0}\tilde{Z}_{j,0}|)^{\alpha/\psi^{ij}}]},
%\end{eqnarray*}
\item
For $d\geq 1$ and $\mathbf{i},\mathbf{j}\in \bbz^d$, the $d$-variate \seq\ $((X_{i_kt}X_{j_kt})_{1\leq k \leq d})_{t \in \mathbb{Z}}$ 
is jointly regularly varying with index $\alpha/\psi^{\mathbf{i},\mathbf{j}}$.
\end{enumerate}
\epr
\bre\label{rem:77} \begin{enumerate}
      \item Equation \eqref{Eq:theta:law:1} shows that in this case the distribution 
of $(\Theta^{i}_t)_{t \geq 0}$ is a mixture of $|\Lambda_i^{(0)}|$ distributions, where 
each distribution gets the weight $1/|\Lambda_i^{(0)}|$. 
Heuristically speaking, a distribution in this mixture that corresponds 
to a specific $(u,v) \in \Lambda_i^{(0)}$ has interpretation as the distribution of $(X_{it}/|X_{i0}|)_{t \geq 0}$, given that we have seen an extreme observation of $|X_{i0}|$ caused by an extreme realization of $\ex^{\eta_{u,v}}$.
The variables $\ex^{\eta_{u,v}}$, $(u,v) \in \Lambda_i^{(0)}$, are those which have a maximum exponent (equal to 1) 
in the product $\prod_{(u,v)}\exp(\psi_{i-u,-v}\eta_{u,v})=\sigma_{i0}$. They are therefore the factors which are most likely 
to make $\sigma_{i0}$, hence $X_{i0}$, extreme. 
\par      
An analogous interpretation can be derived from \eqref{Eq:theta:law:2} for the distribution of $(\Theta^{ij}_t)_{t \geq 0}$. 
      \item Note that for fixed $i, j$, the inner indicator functions in \eqref{Eq:theta:law:1} and \eqref{Eq:theta:law:2} 
are positive only for finitely many $t$. Hence there are only finitely many $t\geq 1$ such that 
$\P(\Theta_t^{i}\neq 0)>0$ and $\P(\Theta_t^{(ij)}\neq 0)>0$.
      \item Using similar techniques as in the proof of cases (1) and (2) below, one can
also give an explicit expression for the resulting $d$-dimensional spectral tail process 
of $((X_{i_kt}X_{j_kt})_{1\leq k \leq d})_{t \in \mathbb{Z}}$ in (3). 
However, due to its complexity, we refrain from stating it here.
     \end{enumerate}
\ere
\begin{proof} 
We start by showing that all mentioned sequences 
are regularly varying. Exemplarily, we show this for case (2). 
Very similar arguments can be used for the two other cases.
For $n \geq 0$ write 
$$ \left(X_{it}X_{jt}\right)'_{t=0, \ldots, n}=\mbox{diag}\left((Z_{it}Z_{jt})_{t=0, \ldots, n}'\right)\cdot
\left(\sigma_{it}\sigma_{jt}\right)'_{t=0, \ldots, n}\,.$$
Since $\psi^{ij} \geq 1$ our moment assumption on $Z$ implies that 
$\E[|Z|^{\alpha/\psi^{ij}+\delta}]<\infty$ for some $\delta>0$. Then  Proposition \ref{prop:rv:1} 
allows us to apply the aforementioned multivariate Breiman lemma, yielding  the  
regular variation of the vector $(X_{it}X_{jt})_{t=0, \ldots, n}$ with index $\alpha/\psi^{ij}$. From the first definition given in Section \ref{subsec:rvts}, this implies the regular variation of the sequence.
\par
As for the derivation of the explicit form of the spectral tail process in (1) and (2), we restrict ourselves to derive the distribution of the spectral tail process $(\Theta_t^{ij})_{t \geq 0}$ in part (2); part (1) is similar.
\par
If $\mu_n^{{\boldsymbol{\sigma}}^{ij}}$ denotes the vague limit measure of 
$(\sigma_{i,0}\sigma_{j,0}, \ldots, \sigma_{i,n}\sigma_{j,n})'$ the multivariate Breiman lemma yields 
the vague limit measure $\mu_n^{\bfX^{ij}}$ of $(X_{i,0}X_{j,0}, \ldots, X_{i,n}X_{j,n})'$ given by
\begin{eqnarray}\nonumber \mu_n^{\bfX^{ij}}(B)&=& c\,\E\left[\mu_n^{\boldsymbol{\sigma}^{ij}}(\times_{t=0}^n (B_t/(Z_{it}Z_{jt})))\right] \\
&=& \label{Eq:mu:X} c\,\E\left[\tilde{\mu}_n^{\boldsymbol{\sigma}^{ij}}\left(\times_{t=0}^n \left(B_t \bigg/ \left(Z_{it}Z_{jt}\prod_{(u,v) \in \Lambda_{i,j,n}^c}\ex^{\eta_{u,v}(\psi_{i-u,t-v}+\psi_{j-u,t-v})}\right)\right)\right)\right] 
\end{eqnarray}
for any $\mu_n^{\bfX^{ij}}$-continuity  Borel set $B=\times_{t=0}^n B_t \in [-\infty,\infty]^{n+1} \setminus \{\mathbf{0}\}$ bounded away from $\mathbf{0}$, $\Lambda_{i,j,n}$ is equal to $\Lambda_{\mathbf{i}, \mathbf{j}, \mathbf{t}}$ 
as defined in \eqref{Eq:Lambda:set:def} with $\mathbf{i}=(i, \ldots, i), \mathbf{j}=(j, \ldots, j), \mathbf{t}=(0, \ldots ,n)$, and $\tilde{\mu}^{\boldsymbol{\sigma}^{ij}}_n$ 
is the limit measure of the \regvary\ vector
\beam \label{Eq:prod:sigma:tilde}
\Big(\prod_{(u,v) \in \Lambda_{i,j,n}}\ex^{\eta_{u,v}\,(\psi_{i-u,t-v}+\psi_{j-u,t-v} )}\Big)_{t=0, \ldots, n}\,,
\eeam
 see the proof of 
Proposition~\ref{prop:rv:1}. The distribution of the tail process of $(X_{it}X_{jt})$ (cf.\ Section \ref{subsec:rvts}) is then determined by
\begin{eqnarray}\label{Eq:theta:law} \P((Y \Theta^{ij}_t)_{t=0, \ldots, n} \in B)&=&\lim_{x \to \infty} \frac{\P((X_{it}X_{jt}/x)_{t=0, \ldots, n} \in B, |X_{i0}X_{j0}|/x>1)}{\P(|X_{i0}X_{j0}|/x>1)}\\
\nonumber &=&\frac{\mu_n^{\bfX^{ij}}(B \cap \big([-\infty,\infty]\backslash [-1,1]\times 
[-\infty,\infty]^n)\big)}{\mu_n^{\bfX^{ij}}\big([-\infty,\infty]\backslash [-1,1]\times [-\infty,\infty]^n\big)}\,.
\end{eqnarray}
The concrete forms of $\tilde{\mu}_n^{\boldsymbol{\sigma}^{ij}}$, hence of $\mu_n^{\mathbf{X}^{ij}}$, 
now depend on whether \eqref{Eq:inf:mom:eta} or \eqref{Eq:conv:eq:eta} holds.
\par
We first assume \eqref{Eq:inf:mom:eta}. Note that
$\Lambda_{i,j,n}=\cup_{t=0}^n \Lambda_{i,j}^{(t)},$ where $\Lambda_{i,j}^{(t)}=\{(u,v): \psi_{i-u,t-v}+\psi_{j-u,t-v}=\psi^{ij}\}$. Indeed, 
we easily see that $\Lambda_{i,j}^{(t)}=\Lambda_{i,j}^{(0)}+(0,t)$, $t=1, \ldots, n$. 
We apply Proposition~\ref{Pr:genRVforproducts}(ii) 
to derive the specific form of the limit measure $\tilde{\mu}_n^{\boldsymbol{\sigma}^{ij}}$ of \eqref{Eq:prod:sigma:tilde}. 
Each component of this vector contains $|\Lambda_{i,j}^{(0)}|$ factors 
with maximal exponent $\psi^{ij}$. For the $t$-th component, those are the factors 
$\exp(\eta_{u,v}(\psi_{i-u,t-v}+\psi_{j-u,t-v}))$, $(u,v) \in \Lambda_{ij}^{(t)}$. Hence $p_{\mbox{\scriptsize eff}}=|\Lambda_{i,j}^{(0)}|$ and 
$P_{\mbox{\scriptsize eff}}=\{\Lambda_{i,j}^{(0)}+(0,t), t=0, \ldots, n\}$. By \eqref{Eq:RV:prod:2}, 
the measure $\tilde{\mu}^{\boldsymbol{\sigma}^{ij}}_n$, up to a constant multiple, is given by
\begin{eqnarray*}
 \tilde{\mu}^{\boldsymbol{\sigma}^{ij}}_n(B)&=& c \sum_{s=0}^n \int_0^\infty \P\bigg(\bigg(\1(\psi_{i-u,t-v}+\psi_{j-u,t-v}=\psi^{ij} \, \forall\, (u,v) \in \Lambda_{i,j}^{(s)})z^{\psi^{ij}}\\
 && \hspace{1.7cm} \prod_{(u,v) \in \Lambda_{i,j,n} \setminus \Lambda_{i,j}^{(s)}}\ex^{\eta_{u,v}(\psi_{i-u,t-v}+\psi_{j-u,t-v})}\bigg)_{0 \leq t \leq n} \in B\bigg)\nu_{\alpha}(dz)\\
 &=& c \sum_{s=0}^n \int_0^\infty \P\bigg(\bigg(\1(t=s)z^{\psi^{ij}}\prod_{(u,v) \in \Lambda_{i,j,n} \setminus \Lambda_{i,j}^{(s)}}\ex^{\eta_{u,v}(\psi_{i-u,t-v}+\psi_{j-u,t-v})}\bigg)_{0 \leq t \leq n} \in B\bigg)\nu_{\alpha}(dz),
\end{eqnarray*}
where $\nu_\alpha(dx)=\alpha x^{-\alpha-1}dx$. 
The $s$-th measure in the sum above is concentrated on the $s$-th axis. Therefore the limit measure $\tilde{\mu}^{\boldsymbol{\sigma}^{ij}}_n$ is
concentrated on the axes.
By \eqref{Eq:mu:X}, this implies that $\mu_n^{\bfX^{ij}}$ is concentrated on the axes as well. Therefore
$\mu_n^{\bfX^{ij}}(B \cap ([-\infty,\infty]\backslash [-1,1])\times [-\infty,\infty]^n)=0$ as soon as one $B_i, 1 \leq i \leq n,$ in $B=\times_{i=0}^n B_i$ is bounded away from 0. 
With \eqref{Eq:theta:law} this gives $Y \Theta_t^{ij}=0$ a.s. for $t \geq 1$ and therefore $\Theta_t^{ij}=0$ a.s. for $t \geq 1$. 
The law of $\Theta_0^{ij}$ follows from the univariate Breiman lemma. 
\par
Next assume \eqref{Eq:conv:eq:eta}. By Proposition~\ref{Pr:genRVforproducts}(i),  
the vague limit measure $\tilde{\mu}^{\boldsymbol{\sigma}^{ij}}_n$ is up to a constant given by 
\begin{eqnarray*} && \tilde{\mu}^{\boldsymbol{\sigma}^{ij}}_n(B) \\
&=&\hspace{-0.3cm} \sum_{(u,v) \in \Lambda_{i,j,n}}\int_0^\infty \P\bigg(\bigg(\1((u,v) \in \Lambda_{i,j}^{(t)})z^{\psi^{ij}}\prod_{\substack{(\tilde{u},\tilde{v}) \in \Lambda_{i,j,n} \\ (\tilde{u},\tilde{v}) \neq (u,v)}}\ex^{(\psi_{i-\tilde{u},t-\tilde{v}}+\psi_{j-\tilde{u},t-\tilde{v}})\eta_{\tilde{u},\tilde{v}}}\bigg)_{t=0, \ldots ,n} \in B \bigg)\nu_{\alpha}(dz).\end{eqnarray*}                                                                                                
For sets $B$ such that $B \cap (\{0\} \times [-\infty,\infty]^n)=\emptyset$ it suffices thereby to sum only 
over  $(u,v) \in \Lambda_{i,j}^{(0)}$ instead over all $(u,v) \in \Lambda_{i,j,n}=\cup_{t=0}^n\Lambda_{i,j}^{(t)}$.
For these sets we have by Breiman's lemma (cf.\ \eqref{Eq:mu:X}),
\begin{align*} & \mu^{\bfX^{ij}}_n(B)/c\\
= &\sum_{(u,v) \in \Lambda_{i,j}^{(0)}}\int\limits_0^\infty \P\bigg((\1((u,v) \in \Lambda_{i,j}^{(t)})z^{\psi^{ij}}\prod\limits_{(\tilde{u},\tilde{v}) \neq (u,v)}\ex^{(\psi_{i-\tilde{u},t-\tilde{v}}+\psi_{j-\tilde{u},t-\tilde{v}})\eta_{\tilde{u},\tilde{v}}} Z_{it}Z_{jt})_{t=0, \ldots ,n} \in B \bigg)\nu_{\alpha}(dz)\\
= &\sum_{(u,v) \in \Lambda_{i,j}^{(0)}}\int\limits_0^\infty \P\bigg((\1((u,v) \in \Lambda_{i,j}^{(t)})z^{\psi^{ij}}X_{it}X_{jt}\ex^{-\psi^{ij}\eta_{u,v}})_{t=0, \ldots ,n} \in B \bigg)\nu_{\alpha}(dz),\end{align*}
where we used that if $(u,v) \in \Lambda_{i,j}^{(t)},$ then 
$$ \prod_{(\tilde{u},\tilde{v}) \neq (u,v)}\ex^{(\psi_{i-\tilde{u},t-\tilde{v}}+\psi_{j-\tilde{u},t-\tilde{v}})\eta_{\tilde{u},\tilde{v}}}= \frac{\sigma_{it}\sigma_{jt}}{\ex^{(\psi_{i-u,t-v}+\psi_{j-u,t-v})\eta_{u,v}}}=\frac{\sigma_{it}\sigma_{jt}}{\ex^{\psi^{ij}\eta_{u,v}}}. $$  
Fubini's Theorem and a substitution finally simplify this expression to 
\begin{eqnarray*} &&  \sum_{(u,v) \in \Lambda_{i,j}^{(0)}}\E\left[\int_0^\infty \1\left(\left(\1((u,v) \in \Lambda_{i,j}^{(t)})z^{\psi^{ij}}X_{it}X_{jt}\ex^{-\psi^{ij}\eta_{u,v}}\right)_{t=0, \ldots, n} \in B\right)\nu_\alpha(dz)\right]  \\
&=& \sum_{(u,v) \in \Lambda_{i,j}^{(0)}}\E\left[\int_0^\infty \1\left(\left(\1((u,v) \in \Lambda_{i,j}^{(t)})y \frac{X_{it}X_{jt}}{|X_{i0}X_{j0}|}\right)_{t=0, \ldots, n} \in B\right)|X_{i0}X_{j0}|^{\alpha/\psi^{ij}}\ex^{-\alpha \eta_{u,v}}\nu_{\frac{\alpha}{\psi^{ij}}}(dy)\right].
\end{eqnarray*}
Note that the range of the inner integral in the last expression can be changed 
from $(0, \infty)$ to $(1,\infty)$, if $B \cap [-1,1] \times [-\infty,\infty]^n = \emptyset$. 
Therefore, by writing
\beao
\tilde{B}_0=B_0 \setminus [-1,1] \,, \quad \tilde{B}_t=B_t, \quad t \geq 1,\quad \tilde B=
\times_{t=0}^n\tilde B_t\,,
\eeao 
we get from \eqref{Eq:theta:law} that
\begin{eqnarray*}
\lefteqn{P\big((Y \Theta^{ij}_t)_{t=0,\ldots,n} \in B\big)}\\&=&\frac{\mu^{\bfX^{ij}}_n(\tilde{B})}{\mu^{\bfX^{ij}}_n(
([-\infty,\infty]\backslash [-1,1])\times [-\infty,\infty]^n)}\\
&=& \frac{\sum\limits_{(u,v) \in \Lambda_{i,j}^{(0)}}\E\left[\int_1^\infty \1\left(\left(\1((u,v) \in \Lambda_{i,j}^{(t)})y \frac{X_{it}X_{jt}}{|X_{i0}X_{j0}|}\right)_{t=0, \ldots, n} \in B\right)|X_{i0}X_{j0}|^{\alpha/\psi^{ij}}\ex^{-\alpha \eta_{u,v}}\nu_{\frac{\alpha}{\psi^{ij}}}(dy)\right]}{\sum\limits_{(u,v) \in \Lambda_{i,j}^{(0)}}\E\left[|X_{i0}X_{j0}|^{\alpha/\psi^{ij}}\ex^{-\alpha \eta_{u,v}}\right]} \\
&=& \sum\limits_{(u,v) \in \Lambda_{i,j}^{(0)}}\frac{1}{| \Lambda_{i,j}^{(0)}|}\frac{\E\left[\1\left(\left(\1((u,v) \in \Lambda_{i,j}^{(t)})Y \frac{X_{it}X_{jt}}{|X_{i0}X_{j0}|}\right)_{t=0, \ldots, n} \in B\right)|X_{i0}X_{j0}|^{\alpha/\psi^{ij}}\right]}{\E\left[|X_{i0}X_{j0}|^{\alpha/\psi^{ij}}\right]},
\end{eqnarray*}
where $Y$ is a Pareto($\alpha/\psi^{ij}$) random variable, independent of all other random variables in the expression. For the last equation, we expanded both numerator and denominator by multiplying with $\E(e^{\alpha \eta_{u,v}})$, noting that for $(u,v) \in \Lambda_{i,j}^{(0)}$ the random variable $e^{\alpha \eta_{u,v}}$ is independent both of the indicator function and of $|X_{i0}X_{j0}|^{\alpha/\psi^{ij}}e^{-\alpha \eta_{u,v}}$. From the law of the tail process $(Y \Theta_t^{ij})$ we can now see that the law of the spectral tail process $(\Theta_t^{ij})$ satisfies \eqref{Eq:theta:law:2}.
\end{proof}

\subsection{Infinite variance stable limit theory for the \sv\ model and its product processes}\label{subsec:limit2}
In the following result we provide central limit theory with infinite variance stable limits for the sums $S_{ij}$; see \eqref{eq:esses}.
\bth\label{thm:1}
We consider the \sv\ model \eqref{eq:1a} and assume the special form of $(\sigma_{it})$ given in \eqref{eq:sv2model} with $\psi=1$. 
For given $(i,j)$, define a \seq\ $(b_n)$ \st\ 
$n\,\P(|X_{i0}X_{j0}|>b_n)\to 1$ as $\nto$.
Assume the following conditions:
\begin{enumerate}
\item 
The conditions of Proposition~\ref{prop:rv:2} hold, ensuring that $\E[|Z|^{\alpha/\psi^{ij}+\vep}]<\infty$ for some $\vep>0$ and
$(X_{it}X_{jt})$ is \regvary\ with index $\alpha/\psi^{ij}$ and spectral tail process $(\Theta_h^{ij})$.
%\item 
%$\E[|Z|^{\alpha/\psi^{ij}+\vep}]<\infty$ for some $\vep>0$ {\red(can be deleted since $\E[|Z|^{\alpha+\vep}]<\infty$ is already assumed in Proposition~\ref{prop:rv:2})},
\item
$(\sigma_{it}\sigma_{jt})$ is $\alpha$-mixing with rate \fct\ $(\alpha_h)$ and there exists $\delta>0$ such that $\alpha_n=o(n^{-\delta})$.
\item Either 
\item[(i)]$\alpha/\psi^{ij}< 1$, or 
\item[(ii)] 
$i\ne j$, $\alpha/\psi^{ij}\in [1,2)$ and $Z$ is symmetric, or 
\item[(iii)] $i=j$,
$\alpha/\psi^{ii}=\alpha/2 \in (1,2)$ and the mixing rate in (2)
satisfies 
$\sup _n n\,\sum_{h=r_n}^\infty \alpha_h<\infty$ for some integer \seq\ $(r_n)$ \st\ $nr_n/b_n^2\to 0$ as $\nto$.
\end{enumerate}
Then 
\beam\label{eq:rra}
b_n^{-1}(S_{ij}-c_n) \std \xi_{ij,\alpha/\psi^{ij}}\,,
\eeam
where $\xi_{ij,\alpha/\psi^{ij}}$ is a totally skewed to the right $\alpha/\psi^{ij}$-stable \rv\ and 
\beao
c_n=\left\{\barr{ll}n\, \E [X^2]& \mbox{$i=j$ and $\alpha\in (2,4)$}\,,\\
0& i\ne j \mbox{ or } \alpha/\psi^{ij}<1 \,,\earr\right.
\eeao
\ethe
\bre\label{reM:c1}
\begin{enumerate}
%Using Breiman's result and the moment condition $\E [|Z|^{\alpha/\psi^{ij}+\vep}]<\infty$ for some $\vep>0$,
%\regvar\ of $(\sigma_{it}\sigma_{jt})$ with index $\alpha/\psi^{ij}$ implies \regvar\ of $(X_{it}X_{jt})$ with the same index. 
%Sufficient conditions for \regvar\ of $(X_{it}X_{jt})$ with index $\alpha/\psi^{ij}$ are given in Proposition~\ref{prop:rv:2}.
\item
If $(\alpha_h)$ decays  at an exponential rate one can choose $r_n= C\log n$ for a sufficiently large constant $C$. Then
  $\sup _n n\,\sum_{h=r_n}^\infty \alpha_h<\infty$ and $nr_n/b_n^2\to 0$ hold. These conditions are also satisfied if
$\alpha_h\le c n^{-(1+\gamma)}$ for some $\gamma>0$, $r_n=C n^\xi$  for some
$\xi>0$ and $1/\gamma\le \xi< 2\psi^{ij}/\alpha-1$.
\item
The \seq\ $(X_{it}X_{jt})$ inherits $\alpha$-mixing from $(\sigma_{it}\sigma_{jt})$; see Remark~\ref{rem:mix}.
\item It is possible to prove joint \con\ for $1\le i,j\le p$ in \eqref{eq:rra}. Due to different tail behavior for
distinct $(i,j)$ the normalizing \seq s $(b_n)=(b_n^{ij})$ typically increase to infinity at different rates. Then it is only of interest
to consider the joint \con\ of those $S_{ij}$ whose summands $X_{it}X_{jt}$ have the same tail index $\alpha/\psi^{ ij}$. More
precisely, it suffices to consider those $S_{ij}$ with the property that $X_{it}X_{jt}$ is tail-equivalent to $X_{it}^2$. 
The joint \con\ follows in a similar way as in the proof below, by observing
that Theorem~\ref{thm:mikwin} is a multivariate limit result. 
The joint limit of $S_{ij}$ in \eqref{eq:rra} with equivalent tails of index $\tilde{\alpha}$ (say) is
jointly $\tilde{\alpha}$-stable with possible dependencies in the limit vector. 
\item
The strongest normalization is needed for $S_i=S_{ii}$. Recall that the summands $X_{it}^2$ of $S_i$ are \regvary\ with index 
$\alpha/2$, i.e., $\psi^{ii}=2$. Let $(a_n)$ be \st\ $n\,\P(|X|>a_n)\to 1$.
Under the conditions of Theorem~\ref{thm:1}, we have that $a_n^{-2}(S_i-c_n)\std \xi_{i,\alpha/2}$, $i=1,\ldots,p$ 
for a jointly $\alpha/2$-stable limit.  If
$\alpha/2<\alpha/\psi^{ij}$ for some $i\ne j$, then $b_n/a_n^2\to 0$, hence $a_n^{-2} S_{ij}\stp 0$. It is possible that $X_{it}X_{jt}$
is regularly varying with index $\alpha/2$ but nevertheless $b_n/a_n^2\to 0$; see Example~\ref{exam:1} 
which deals with the case $\E [\ex^{\alpha \eta}]=\infty$.
\end{enumerate}
\ere
\begin{proof}
We apply Theorem~\ref{thm:mikwin} to the \seq\ $(X_{it}X_{jt})$, cf.\ also Remark~\ref{rem:case:alpha:1}.\\[2mm] 
(1) The \regvar\ condition on $(X_{it}X_{jt})$ with index $\alpha/\psi^{ij}$ is satisfied by assumption. Moreover, $\Theta_h=0$ for 
sufficiently large $h$; see Remark~\ref{rem:77}.\\[2mm]
(2) The assumption about the mixing coefficients in condition (2) implies that for a 
sufficiently small $\varepsilon \in (0,1)$ and $m_n=n^{1-\varepsilon}$ there exists an integer sequence $l_n=o(m_n)$ such that $k_n \alpha_{l_n} \to 0$. For this choice of $m_n$ and $l_n$, the proof of the mixing condition for the sums of the truncated variables
\beao
\un S_{ij} =\sum_{t=1}^n X_{it}X_{jt} \1( |X_{it}X_{jt}|>\vep b_n)
\eeao
is now analogous to the proof of the corresponding property in Theorem~\ref{the:1}. \\[2mm]
(3) We want to show that
\beam\label{eq:ac}
\lim_{l\to\infty}\limsup_{\nto} n\,\sum_{t=l}^{m_n}\,\P\big(|X_{it}X_{jt}|>b_n\,,|X_{i0}X_{j0}|>b_n\big)=0\,
\eeam
for $m_n=n^{1-\varepsilon}$ as above.
Write
$$ 
\sigma_{it}\sigma_{jt}= \prod_{(m,n)}\exp((\psi_{i-m,t-n}+\psi_{j-m,t-n})\eta_{m,n})
$$
and set $\Lambda_{\varepsilon,t}=\{(m,n): \psi_{i-m,t-n}+\psi_{j-m,t-n}\geq 8^{-1}\psi^{ij}\varepsilon\}$, $t \in \mathbb{Z}$. 
Without loss of generality we assume that $l$ is so large that $\Lambda_{\varepsilon,t} \cap \Lambda_{\varepsilon,0}$ is empty
for all $t \geq l$. Then write for $t \geq l$, 
\beao
\sigma_{it}\sigma_{jt}=
\sigma_{it,jt,\Lambda_{\varepsilon,t}}\cdot \sigma_{it,jt,\Lambda_{\varepsilon,0}} \cdot \sigma_{it,jt,\Lambda_{\varepsilon,0,t}^c},\;\;\; \sigma_{i0}\sigma_{j0}=\sigma_{i0,j0,\Lambda_{\varepsilon,0}}\cdot \sigma_{i0,j0,\Lambda_{\varepsilon,t}} \cdot \sigma_{i0,j0,\Lambda_{\varepsilon,0,t}^c},\eeao
where
$$\sigma_{it_1,jt_1,\Lambda_{\varepsilon,t_2}}= \prod_{(m,n) \in \Lambda_{\varepsilon,t_2}}\exp((\psi_{i-m,t_1-n}+\psi_{j-m,t_1-n})\eta_{m,n}).$$
We conclude that $(\sigma_{it,jt,\Lambda_{\varepsilon,t}}$,$\sigma_{it,jt,\Lambda_{\varepsilon,0}}$, $\sigma_{i0,j0,\Lambda_{\varepsilon,0}}$, $\sigma_{i0,j0,\Lambda_{\varepsilon,t}})$ and $(\sigma_{it,jt,\Lambda_{\varepsilon,0,t}^c}$,$\sigma_{i0,j0,\Lambda_{\varepsilon,0,t}^c})$
are independent. We have
\beao\lefteqn{
\P\big(|X_{it}X_{jt}|>b_n\,,|X_{i0}X_{j0}|>b_n\big)}\\
&\le & \P\big(\max(|Z_{i0}Z_{j0}|,|Z_{it}Z_{jt}|)\,\max(\sigma_{it,jt,\Lambda_{\varepsilon,0,t}^c},\sigma_{i0,j0,\Lambda_{\varepsilon,0,t}^c})\, \\&&
\min(\sigma_{i0,j0,\Lambda_{\varepsilon,0}} \sigma_{i0,j0,\Lambda_{\varepsilon,t}},\sigma_{it,jt,\Lambda_{\varepsilon,t}}\sigma_{it,jt,\Lambda_{\varepsilon,0}})>b_n\big)\,.
%&\le &
%\P\big(M^2\max(\sigma_{i0,\Lambda_2^c}\sigma_{j0,\Lambda_2^c},\sigma_{it,\Lambda_2^c}\sigma_{jt,\Lambda_2^c})\, \min(\sigma_{i0,\Lambda}\sigma_{j0,\Lambda},\sigma_{it,\Lambda}\sigma_{jt,\Lambda})>b_n\big)\\
%&&+\P\big(\{\max(\sigma_{i0,\Lambda^c}\sigma_{j0,\Lambda^c},\sigma_{it,\Lambda^c}\sigma_{jt,\Lambda^c})\, \min(\sigma_{i0,\Lambda}\sigma_{j0,\Lambda},\sigma_{it,\Lambda}\sigma_{jt,\Lambda})>b_n\}\\&&\cap 
%\big(\{\sigma_{i0,\Lambda_1^c}\sigma_{j0,\Lambda_1^c}>M\}\cup \{\sigma_{it,\Lambda_1^c}\sigma_{jt,\Lambda_1^c}>M\}\big)=:J_{1t}+
%J_{2t}\,.
\eeao
The \ds\ of $\max(\sigma_{it,jt,\Lambda_{\varepsilon,0,t}^c},\sigma_{i0,j0,\Lambda_{\varepsilon,0,t}^c})$ 
is stochastically dominated uniformly for $t\geq l$ by a \ds\ which has moment of order $8 \alpha/(\psi^{ij}\varepsilon)>2 \alpha/\psi^{ij}$. Furthermore,
\begin{eqnarray*}
 && \min(\sigma_{i0,j0,\Lambda_{\varepsilon,0}} \sigma_{i0,j0,\Lambda_{\varepsilon,t}},\sigma_{it,jt,\Lambda_{\varepsilon,t}}\sigma_{it,jt,\Lambda_{\varepsilon,0}})\\
 &\le & \min\Big(\prod_{(m,n) \in \Lambda_{\varepsilon,0} \cup \Lambda_{\varepsilon,t}}\exp((\psi_{i-m,-n}+\psi_{j-m,-n})(\eta_{m,n})_+),\\
 && \hspace{1cm}\prod_{(m,n) \in \Lambda_{\varepsilon,0} \cup \Lambda_{\varepsilon,t}}\exp((\psi_{i-m,t-n}+\psi_{j-m,t-n})(\eta_{m,n})_+)\Big)\,\\
 &\le& \min\Big(\prod_{(m,n) \in \Lambda_{\varepsilon,0}}\exp(\psi^{ij}(\eta_{m,n})_+)\prod_{(m',n') \in \Lambda_{\varepsilon,t}}\exp(8^{-1}\psi^{ij}\varepsilon (\eta_{m',n'})_+),\\
 && \hspace{1cm}\prod_{(m',n') \in \Lambda_{\varepsilon,t}}\exp(\psi^{ij}(\eta_{m',n'})_+)\prod_{(m,n) \in \Lambda_{\varepsilon,0}}\exp(8^{-1}\psi^{ij}\varepsilon (\eta_{m,n})_+)\Big)\, \\
 &\leq& \min\Big(\prod_{(m,n) \in \Lambda_{\varepsilon,0}}\exp((\psi^{ij}+8^{-1}\psi^{ij}\varepsilon)(\eta_{m,n})_+),
 \prod_{(m,n) \in \Lambda_{\varepsilon,t}}\exp((\psi^{ij}+8^{-1}\psi^{ij}\varepsilon)(\eta_{m,n})_+)\Big).
 \end{eqnarray*}
The \rhs\ is \regvary\ with index $2\alpha/(\psi^{ij}(1+8^{-1}{\vep}))$.
A stochastic domination argument and an application of Breiman's lemma show that uniformly for $l\le t\le m_n$,
$$ m_n\, n\, \P\big(|X_{it}X_{jt}|>b_n\,,|X_{i0}X_{j0}|>b_n\big)=n^{2-\varepsilon} o\left(b_n^{-2\alpha/(\psi^{ij}(1+4^{-1}\epsilon))}\right)=n^{2-\varepsilon}o(n^{-2/(1+2^{-1}\varepsilon)})=o(1)$$
which yields \eqref{eq:ac}.\\[2mm]
(4) We check the vanishing small values condition. For any fixed $\delta$, we write
\beao
\ov {X_{it}X_{jt}}&=& X_{it}X_{jt} \1(|X_{it}X_{jt}|\le \delta b_n)\,,\qquad i\ne j\,,\\
\ov {X_{it}^2}&=& X_{it}^2 \1(X_{it}^2\le \delta b_n) - \E[X_{it}^2 \1(X_{it}^2\le \delta b_n)] \,,\\
\ov S_{ij}&=& \sum_{t=1}^n \ov {X_{it}X_{jt}}\,,\qquad \ov S_i=\ov S_{ii}\,.
\eeao
%Assume $\alpha/\psi^{ij}\in(0,1)$.
%We have by Karamata's theorem for any $\gamma>0$ as $\nto$.
%\beao
%\P(|\ov S_{ij}|>\gamma b_n)&\le &
%(\gamma b_n)^{-1} \E [\ov S_{ij}]\\&\le& n\,(\gamma b_n)^{-1}\,\E[|\ov {X_{it}X_{jt}}|]\\
%&\sim & \gamma^{-1} \delta^{1-\alpha/\psi^{ij}}\,,
%\eeao
%and the \rhs\ converges to zero for $\delta\downarrow 0$. {\red Can be deleted since index is less than 1 and vanishing small values condition does not need to hold.}
Assume $\alpha/\psi^{ij}\in [1,2)$, $i\ne j$. Then, by symmetry of the \rv s $Z_{it}$ and Karamata's theorem for any $\gamma>0$ as $\nto$,
$\E[\ov S_{ij}]=0$ and
\beao
\P(|\ov S_{ij}|>\gamma b_n)&\le &
(\gamma b_n)^{-2} \var(\ov S_{ij})\\
&=& n\,(\gamma b_n)^{-2} \E[(\ov {X_{it}X_{jt}})^2]\\
&\sim &\gamma^{-2}\,\delta^{2-\alpha}\,,
\eeao
and the \rhs\ converges to zero as $\delta\downarrow 0$.
\par
For $i=j$ and $\alpha/\psi^{ii}>1$ we need a different argument. 
We have by \v Cebyshev's inequality,
\beao
\P(|\ov S_{i} |>\gamma \,b_n)&\le & \gamma^{-2} b_n^{-2} \var\big(\ov S_{i}\big)\\
&=&  \gamma^{-2}\,(n/b_n^2)\,\sum_{|h|<n} (1-h/n)\,\cov(\ov {X_{i0}^2}, \ov {X_{ih}^2})\,.
\eeao
For $|h|\le h_0$ for any fixed $h_0$, $(n/b_n^2)  |\cov(\ov {X_{i0}^2}, \ov {X_{ih}^2})|$ vanishes by letting first $\nto$ and then $\delta\downarrow 0$.
This follows by Karamata's theorem. Standard bounds for the covariance \fct\ of an $\alpha$-mixing \seq\ 
(see Doukhan \cite{doukhan:1994}, p.~3)
yield
\beao
(n/b_n^2)\sum_{r_n\le |h|<n}|\cov(\ov {X_{i0}^2}, \ov {X_{ih}^2})| 
&\le & c\,\delta^2 n\,\sum_{r_n\le |h|<n}\alpha_h\,,
\eeao
where $r_n\to\infty$ is chosen \st\ $\sup_n\,n\,\sum_{r_n\le |h|<\infty} \alpha_h<\infty$ and $nr_n/b_n^2 \to 0$. The \rhs\ converges to zero by first letting $\nto$ and then
$\delta\downarrow 0$. It remains to show that
\beao
I_n=(n/b_n^2)\sum_{h_0< |h|\le r_n}(1-h/n)\,\cov(\ov {X_{i0}^2}, \ov {X_{ih}^2})
\eeao
is \asy ally negligible. We have
\beao
|I_n|&\le& (n/b_n^2)\sum_{h_0< |h|\le r_n} \E [X_{i0}^2\,X_{ih}^2\1(X_{i0}^2\le \delta b_n,X_{ih}^2\le \delta b_n)] + 
c\, n\,r_n/b_n^2\\
&\le & (n/b_n^2)\sum_{h_0< |h|\le r_n} \E [X_{i0}^2\,X_{ih}^2] +o(1)\,,
\eeao 
where we used that $n\,r_n/b_n^2\to 0$.
We will show that the summands on the \rhs\ are uniformly bounded by a constant if $h_0$ is sufficiently large. Then $\lim_{\nto} I_n=0$. 
\par
We observe that by H\"older's inequality,
\beao
\E [X_{i0}^2\,X_{ih}^2]&=& c\, \E [\sigma_{i0}^2\,\sigma_{ih}^2]\\
&=&c\,\E \big[
\ex^{2\sum_{(k,l)\in \Gamma_{\xi}}\psi_{kl}(\eta_{i-k,-l}+\eta_{i-k,h-l})}  \ex^{2\sum_{(k,l) \not \in \Gamma_\xi}\psi_{kl}(\eta_{i-k,-l}+\eta_{i-k,h-l})}\big]\\
&\le &c\,\big(\E \big[
\ex^{2r\sum_{(k,l)\in \Gamma_\xi}\psi_{kl}(\eta_{i-k,-l}+\eta_{i-k,h-l})}\big]\big)^{1/r} \big(\E\big[ \ex^{2s\sum_{(k,l) \not \in \Gamma_\xi}\psi_{kl}(\eta_{i-k,-l}+\eta_{i-k,h-l})}\big]\big)^{1/s}\,,
\eeao
where $\Gamma_\xi=\{(k,l): \psi_{ik}>\xi\}$ for some positive $\xi,s,t$ \st\ $1/r+1/s=1$. Since $\sigma_{i0}^2$ has moments up to order $\alpha/\psi^{ii}\in (1,2)$ and 
$(\eta_{i-k,-l})_{(k,l)\in\Gamma_\xi}$ and $(\eta_{i-k,h-l})_{(k,l)\in \Gamma_\xi}$ are independent for sufficiently large $h$ we can choose $r>1$ close to
one \st\ $\E \big[
\ex^{2r\sum_{(k,l)\in \Gamma\xi}\psi_{kl}(\eta_{i-k,-l}+\eta_{i-k,h-l})}\big]$ is finite. This implies that we choose $s$ sufficiently large. On the other hand,
for fixed $s$ we can make $\xi$ so small that $\E\big[ \ex^{2s\sum_{(k,l) \not \in \Gamma\xi}\psi_{kl}(\eta_{i-k,-l}+\eta_{i-k,h-l})}\big]$ is finite and uniformly bounded for
sufficiently large $h$. Fine tuning $\xi$ and $s$, we may conclude that $\lim_{\nto}I_n=0$ as desired.

By Theorem~\ref{thm:mikwin} and Remark~\ref{rem:case:alpha:1} the result now follows; see  also 
the end of the proof of Theorem~\ref{the:1} for the form of the resulting limit law.
\end{proof}

\bexam\label{exam:1}
We assume that $\E [\ex^{\alpha \eta}]=\infty$, hence $\ex^{2\eta}$ does not have a finite $\alpha/2$-th moment. 
Using Lemma~\ref{lem:product}(5), calculation shows that for $i\ne j$ with $\psi^{ij}=2$,
\beam\label{eq:null}
\lim_{\xto}\dfrac{\P(|X_{i0}\,X_{j0}|>x)}{\P(X^2>x)}=0
\eeam
Define $(a_n)$ \st\ $n\,\P(|X|>a_n)\to 1$. We may conclude from \eqref{eq:null} and Theorem~\ref{thm:1} that for $i\ne j$ we have 
$a_n^{-2} S_{ij}\stp 0$ although 
both $X_{i0}\,X_{j0}$ and $X^2$ are \regvary\ with index $\alpha/2$.
\par
By Theorem~\ref{thm:1} and Remark~\ref{reM:c1} we conclude that
\beam\label{eq:aa}
a_n^{-2} (S_i-c_n)_{i=1,\ldots,p} \std (\xi_{i,\alpha/2})_{i=1,\ldots,p}\,,
\eeam 
where the limit vector consists of $\alpha/2$-stable components. The spectral tail process $(\bfTh_h)_{h\ge 1}$ 
of the \seq\ $\bfX_t=(X_{1t},\ldots,X_{pt})'$, $t=1,2,\ldots$, vanishes. This follows by an argument similar to the proofs of 
Propositions~\ref{prop:rv:2} and \ref{Pr:genRVforproducts} under condition \eqref{Eq:inf:mom:eta}. 
A similar argument also yields that
\beao
\lim_{\xto}\dfrac{\P(|X_{i0}|>x\,,|X_{j0}|>x)}{\P(|X|>x)}=0\,,\qquad i\ne j\,.
\eeao
Therefore the the \ds\ of $\bfTh_0$ is concentrated on the axes and has the same form as $\boldsymbol{\Theta}_0^{(2)}$ in \eqref{Eq:Theta_0:measure}. As in the proof of Theorem~\ref{the:1} this implies that the limit random vector in \eqref{eq:aa} has iid components.
\par
We conclude that the limit theory for $S_{ij}$, $1\le i,j\le p$, are very essentially the same 
in Case (1) and in Case (2) when the additional
condition $\E [\ex^{\alpha \eta}]=\infty$ holds.
\eexam
\bexam\label{exam:2}
Assume that \eqref{Eq:conv:eq:eta} holds. We may conclude from Theorem~\ref{thm:1}
that $a_n^{-2} S_{ij}\stp 0$ for $i\ne j$ if $\psi^{ij}<2$. The crucial difference to the previous case appears
when $\psi^{ij}=2$ for some $i\ne j$. In this case, not only the $(a_n^{-2}(S_i-c_n))$, $i=1,2,\ldots,$ have 
totally skewed to the right $\alpha/2$-stable limits 
but we also have $a_n^{-2} S_{ij}\std \xi_{ij,\alpha/2}$ for non-degenerate $\alpha/2$-stable $\xi_{ij,\alpha/2}$. From \eqref{eq:psi:ij:def}
we conclude that if $\psi^{ij}=2$ appears then $\psi^{i'j'}=2$ for all $|i'-j'|=|i-j|$. This means that non-degenerate limits may appear not only
on the diagonal of the matrix $a_n^{-2}(S_{ij}-c_n)$ but also along full sub-diagonals.
\par
In this case, the distribution of $\boldsymbol{\Theta}_0$ from the spectral tail process of the 
\seq\ $\bfX_t=(X_{1t},\ldots,X_{pt})'$ does not have to be concentrated on the axes---in contrast to Example~\ref{exam:1}. 
This implies that the limiting $\alpha/2$-stable random variables $\xi_{i,\alpha/2}, i=1, \ldots, p,$ are in general not independent. However, similar to the arguments at the end of the proof of Theorem~\ref{the:1}, one can show that the distribution of the limiting random vector $(\xi_{i,\alpha/2})_{i=1,\ldots,p}$ is the convolution of distributions of $\alpha/2$-stable random vectors which concentrate on 
hyperplanes of $\mathbb{R}^p$ of dimension less or equal than $|\{(m,n): \psi_{mn}=1\}|$. 
\eexam
\subsection{The eigenvalues of the sample covariance matrix of a multivariate \sv\ model}\label{subsec:eigen2}
In this section we provide some results for the eigenvalues of the  sample covariance matrix $\bfX^n(\bfX^n) '$
under the conditions of  Theorem~\ref{thm:1}.
We introduce the sets 
\beao
\Gamma_p= \{(i,j):1\le i,j\le p\;\mbox{\st}\; \psi^{ij}=2\}\,,\qquad \Gamma_p^c= \{(i,j):1\le i,j\le p\}\backslash \Gamma_p\\
\eeao 
and let $(a_n)$ be \st\ $n\,\P(|X|>a_n)\to 1$. 
\bth\label{thm:2}
Assume that the conditions of Theorem~\ref{thm:1} hold for $(X_{it},X_{jt})$, $1\le i,j\le p$, and $\alpha\in (0,4)$.
Then
\beao
a_n^{-2} \big\|\bfX^n(\bfX^n)'- \wt \bfX^n\big\|_2\stp 0\,,\qquad \nto\,,
\eeao
where $\wt \bfX^n$ is a $p\times p$ matrix with entries
\beao
\wt X_{ij}=\sum_{t=1}^n X_{it}X_{jt} \1( (i,j)\in \Gamma_p)\,,\qquad 1\le i,j\le p\,.
\eeao
Moreover, if $\E [\ex^{\alpha\eta}]=\infty$ we also have
\beao
a_n^{-2} \big\|\bfX^n(\bfX^n)'- {\rm diag}(\bfX^n(\bfX^n)')\big\|_2\stp 0\,,\qquad \nto\,.
\eeao
\ethe
\begin{proof} We have
\beao
a_n^{-4} \big\|\bfX^n(\bfX^n)'- \wt \bfX^n\big\|_2^2 \le 
\sum_{(i,j)\in \Gamma_p^c} \big(a_n^{-2}S_{ij}\big)^2 \,.
\eeao
For $(i,j) \in \Gamma_p^c$ we have $i\ne j$ and the \seq\ $(X_{it}X_{jt})$ is \regvary\ with index $\alpha/\psi^{ij}>\alpha/2$. 
In view of Theorem~\ref{thm:1} the \rhs\
converges to zero in \pro y. 
\par
In the case when $\E[\ex^{\alpha\eta}]=\infty$ we learned in Example~\ref{exam:1} that $a_n^{-2}S_{ij}\stp 0$ whenever $i\ne j$. This concludes
the proof.
\end{proof}
For any $p\times p$ non-negative definite matrix $\bfA$ write $\la_i(\bfA)$, $i=1,\ldots,p$, for its eigenvalues and 
$\la_{(1)}(\bfA)\ge \cdots \ge \la_{(p)}(\bfA)$ for their ordered values. For the eigenvalues of $\bfX^n(\bfX^n)'$ we keep the previous
notation $(\la_i)$,
\bco
Assume the conditions of Theorem~\ref{thm:2} and $\alpha\in (0,4)\backslash \{2\}$. Then 
\beam\label{eq:xy}
a_n^{-2}\max_{i=1,\ldots,p} \big|\la_{(i)}-\la_{(i)}(\wt \bfX^{n})\big|\stp 0\,.
\eeam
and
\beam\label{eq:90}
a_n^{-2}\Big(\la_{(i)} -n\,\E[X^2] 1(\alpha\in (2,4))\Big)_{i=1,\ldots,p}
\std \Big(\la_{(i)}\big((\xi_{kl,\alpha/2} \1((k,l)\in \Gamma_p))_{1\le k,l\le p}\big)\Big)_{i=1,\ldots,p} \,,\nonumber\\
\eeam
where $(\xi_{ij,\alpha/2})_{(i,j)\in \Gamma_p}$ are jointly $\alpha/2$-stable (possibly degenerate for $i\ne j$) \rv s.
Moreover, in the case when $\E [\ex^{\alpha\eta}]=\infty$ we have
\beam\label{eq:90a}
a_n^{-2}\Big(\la_{(i)} -n\,\E[X^2] 1(\alpha\in (2,4))\Big)_{i=1,\ldots,p}
\std \big(\xi_{(i),\alpha/2}\big)_{i=1,\ldots,p}\,,
\eeam
where $(\xi_{i,\alpha/2})_{i=1,\ldots,p}$ are iid totally skewed to the right $\alpha/2$-stable \rv s with order statistics $\xi_{(1),\alpha/2}\ge \cdots\ge \xi_{(p),\alpha/2}$.
\eco
\begin{proof} Relation \eqref{eq:xy} is an immediate con\seq\ of Theorem~\ref{thm:2} and Weyl's inequality; 
see Bhatia \cite{bhatia:1997}. We conclude  from Theorem~\ref{thm:1} and Remark~\ref{reM:c1}(3) that
\beam\label{eq:nn}
a_n^{-2}\big(S_{ij}  - n\,\E[X^2]\,\1 (\alpha\in (2,4))\big)_{(i,j)\in \Gamma_p}\std \big(\xi_{ij,\alpha/2}\big)_{(i,j)\in \Gamma_p} \,.
\eeam
Then \eqref{eq:90} follows. Relation \eqref{eq:90a} is a special case of \eqref{eq:90}. If $\E [\ex^{\alpha\eta}]=\infty$ then, 
in view of Example~\ref{exam:1}, only the diagonal elements in \eqref{eq:nn} have non-degenerate iid $\alpha/2$-stable limits.
\end{proof}
\subsubsection*{Some conclusions}
By virtue of this corollary and in view of Section~\ref{subsec:eigen} the results for the eigenvalues 
in Case (1) and in Case (2) when $\E [\ex^{\alpha\eta}]=\infty$ are very much the same. Moreover, the results in 
Section~\ref{subsec:appl} remain valid in the latter case.
\par
If \eqref{Eq:conv:eq:eta} holds, Case (2) is quite different from Case (1); see Example~\ref{exam:2}.
In this case not only the diagonal of the matrix $\bfX^n(\bfX^n)'$ determines the
\asy\ behavior of its eigenvalues and eigenvectors. Indeed, if $\psi^{ij}=2$ for some $i\ne j$, then
at least two sub-diagonals of $\bfX^n(\bfX^n)'$ have non-degenerate $\alpha/2$-limits and these sub-diagonals together
with the diagonal determine the \asy\ behavior of the eigenspectrum. The limiting diagonal elements are dependent
in contrast to Case (1). This fact and the presence of sub-diagonals are challenges if one wants to calculate
the limit \ds s of the eigenvalues and eigenvectors.
\section{Simulations and data example}\label{sec:simulation}
\begin{figure}[h]
\includegraphics[width=8cm]{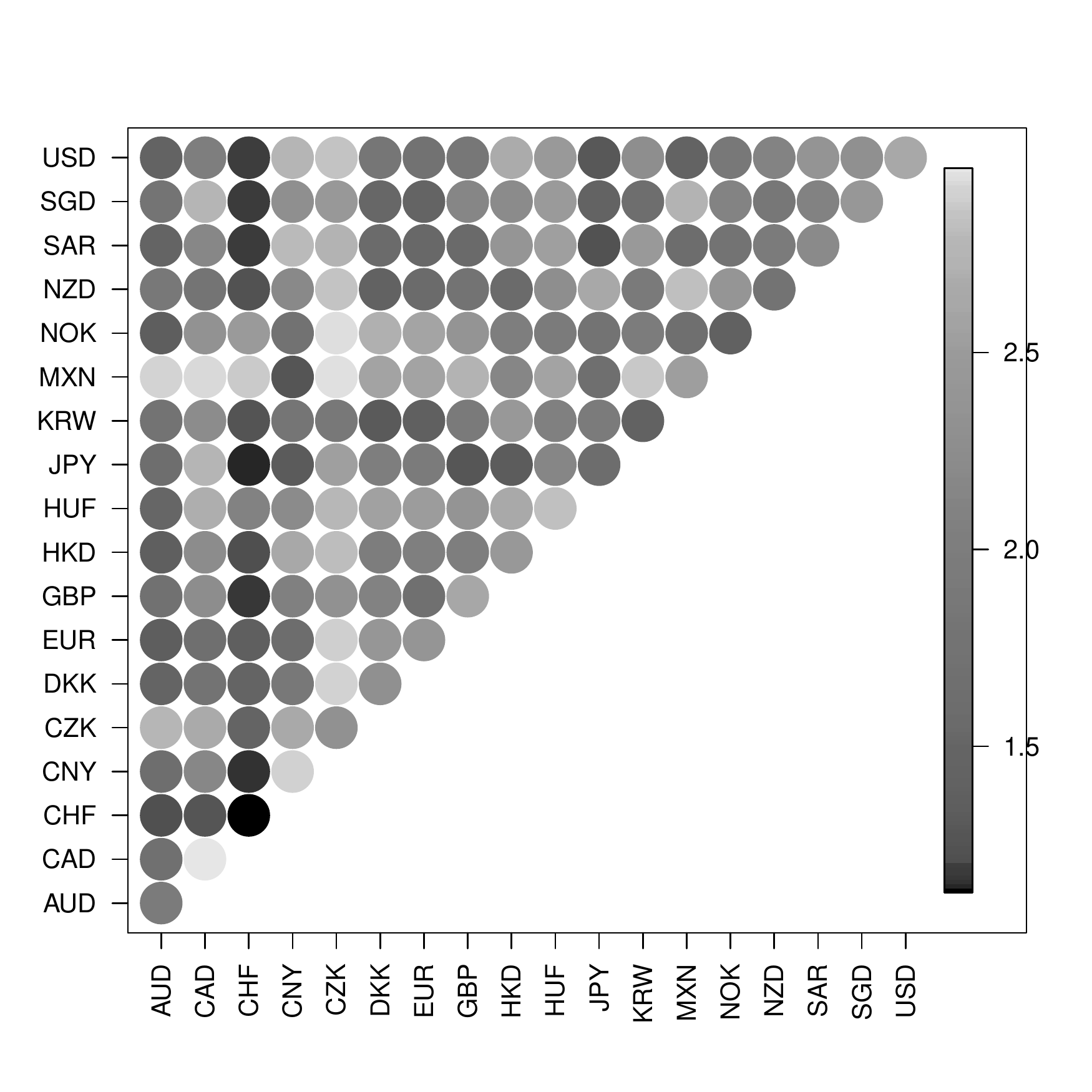}
\caption{Estimated tail indices of cross products for the FX rates of 18 currencies against SEK. The indices are derived by Hill estimators with threshold equal to the 97\%-quantile of $n=1567$ observations.}
\label{Fig:Hill}
\end{figure}

\setcounter{equation}{0}
In this section we illustrate the behavior of sample covariance matrices for moderate sample sizes
for the models discussed in Sections~\ref{sec:case1} and \ref{Sec:case2} 
and we compare them with a real-life data example. These data consist of
1567 daily log-returns of foreign exchange (FX) rates from
18 currencies against the Swedish Kroner (SEK) from January 4th 2010 to April 1st 2016, as made available by the Swedish National Bank. To start with, the Hill estimators of the tail indices $\alpha_{ij}, 1 \leq i,j, \leq 18,$ of the cross 
products $X_{it}X_{jt}, 1 \leq i,j, \leq 18,$ are visualized in Figure \ref{Fig:Hill}. In particular, 
the Hill estimators on the diagonal (corresponding to the series $X_{it}^2, 1 \leq i \leq 18$) of the values $\alpha_i/2$, where $\alpha_i$ is the tail index of the $i$th currency,
are of similar size although not identical. Even if all series had the same tail index the Hill estimator
exhibits high statistical uncertainty which even increases for serially dependent data, cf.\ Drees \cite{drees:2003}.
A way to make the data more homogeneous in their tails is to rank-transform their marginals to the same \ds . We do, however,
refrain from such a  transformation to keep the correlation structure of the original data unchanged. 

It is clearly visible that some off-diagonal components of the matrix have an estimated tail index 
which is comparable to the on-diagonal elements. 
This implies that the tails of the corresponding 
off-diagonal entries $S_{ij}, i\ne j$, of the sample covariance matrix may 
be of a similar magnitude as the on-diagonal entries $S_i$. This is in stark contrast to the 
asymptotic behavior of the models analyzed in Section \ref{sec:case1}. 

\begin{figure}[h]
\subcaptionbox{Based on FX rate data of 18 foreign currencies against SEK.\label{eigen:data}}
{\includegraphics[width=6.5cm]{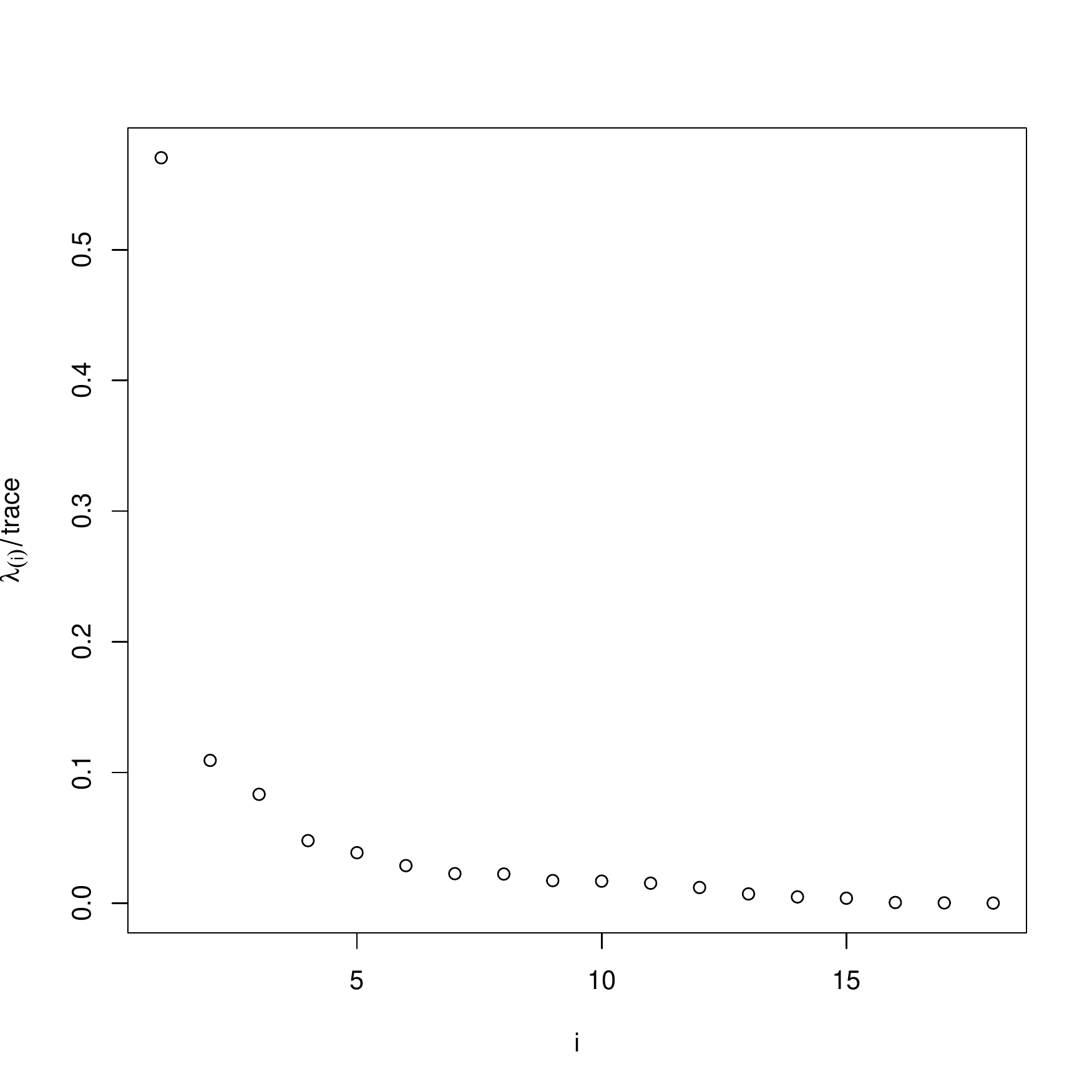}
\includegraphics[width=6.5cm]{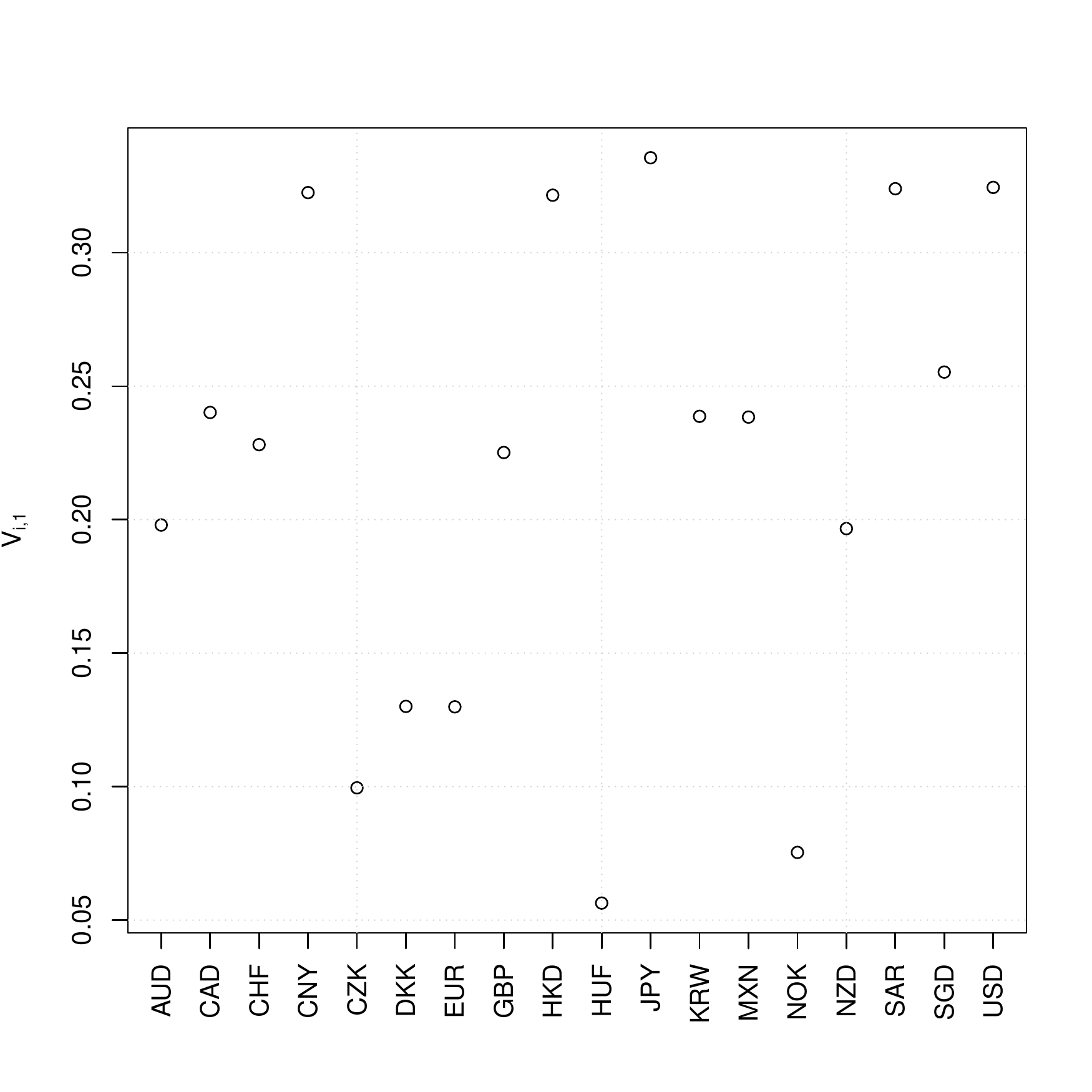}}
\caption{Normalized and ordered eigenvalues (left) and eigenvector corresponding to largest eigenvalue (right) of real and simulated data, with $n=1567, p=18$.}\label{fig:1}
\end{figure}
\begin{figure}
\ContinuedFloat
\subcaptionbox{Based on a \sv\ model with heavy-tailed innovation sequence.\label{eigen:model1}}
{\includegraphics[width=6.5cm]{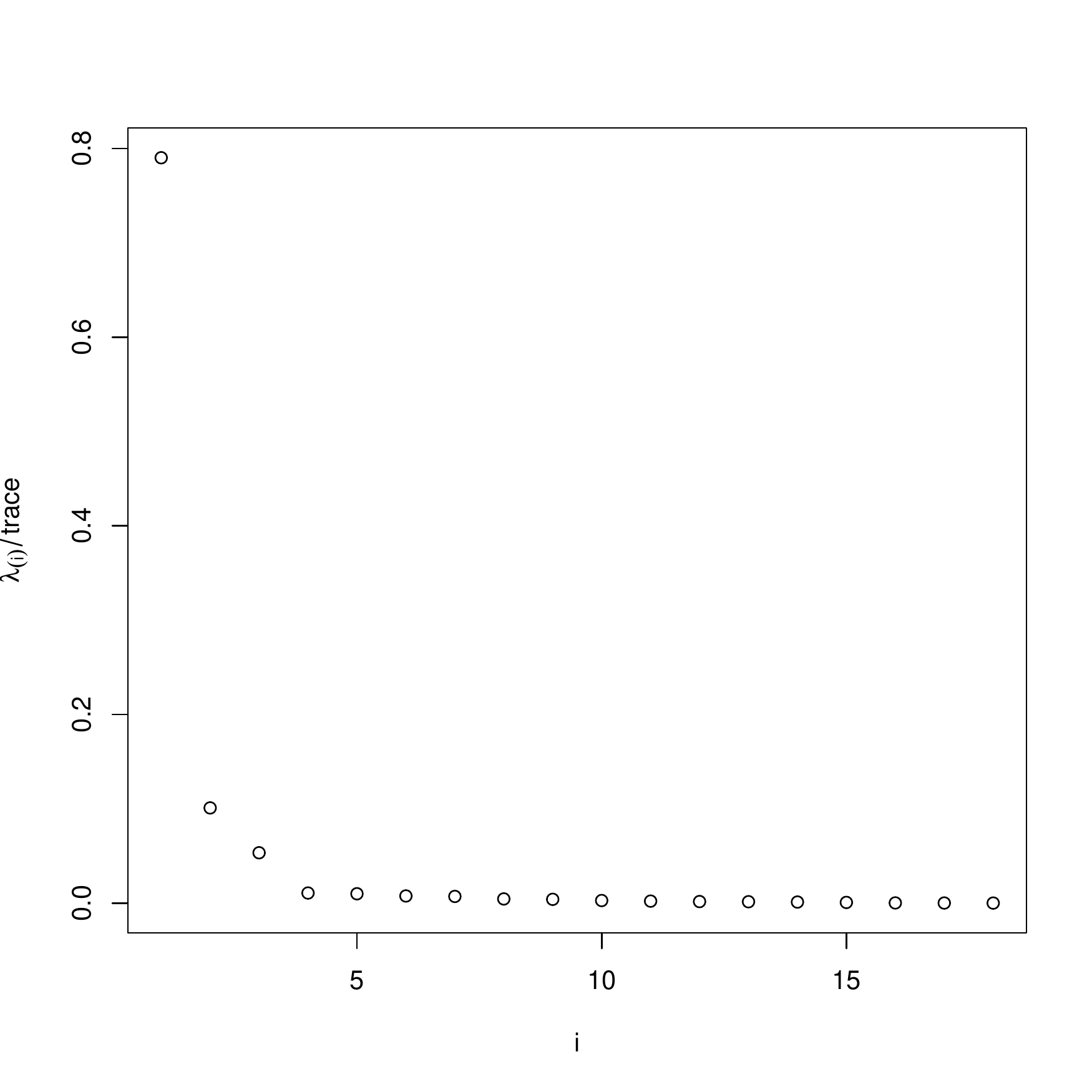}
\includegraphics[width=6.5cm]{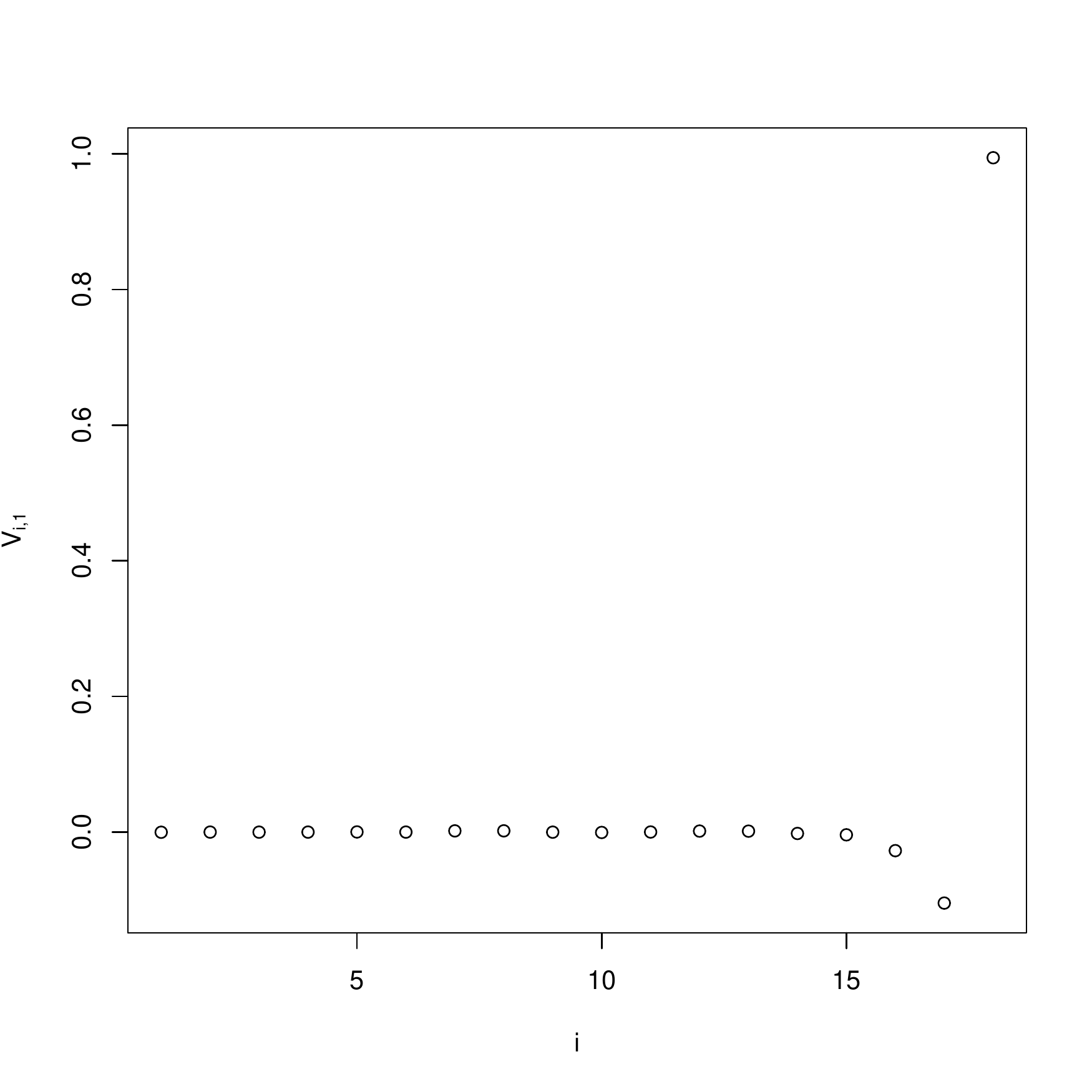}}
\subcaptionbox{Based on a \sv\ model with heavy-tailed volatility sequence that satisfies assumptions of Example \ref{exam:1}. \label{eigen:model2}}
{\includegraphics[width=6.5cm]{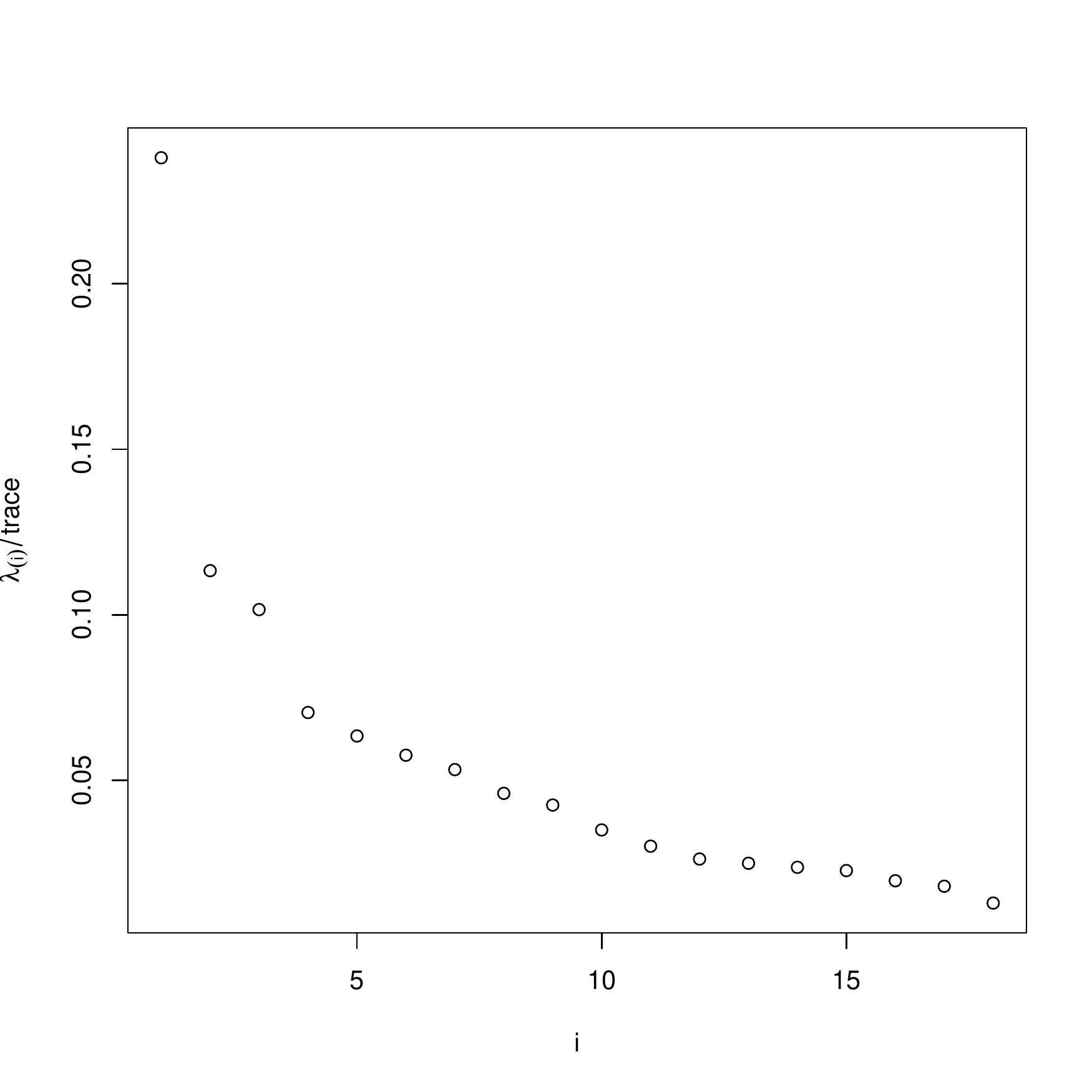}
\includegraphics[width=6.5cm]{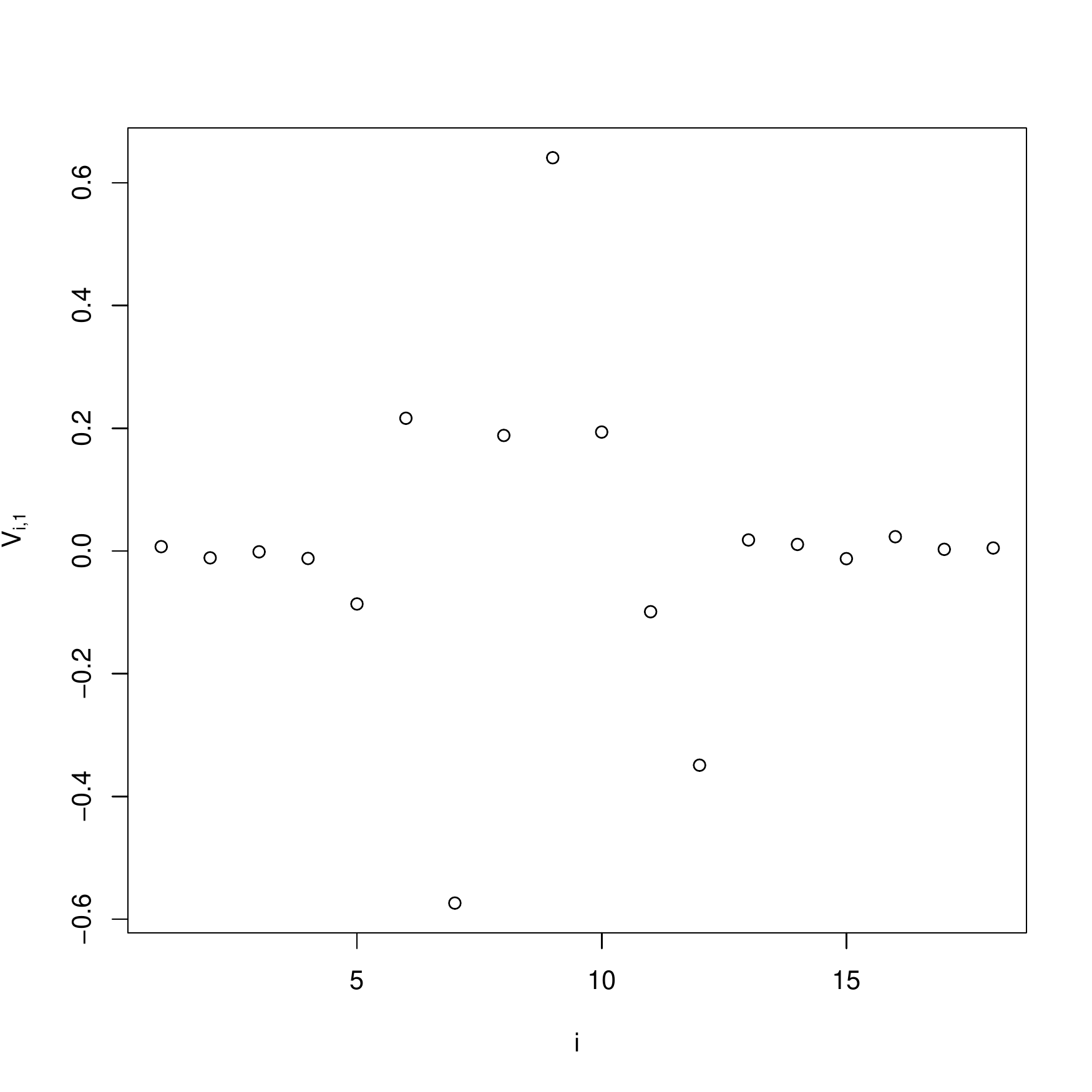}}
\subcaptionbox{Based on a \sv\ model with heavy-tailed volatility sequence that satisfies assumptions of Example \ref{exam:2}.  \label{eigen:model3}}
{\includegraphics[width=6.5cm]{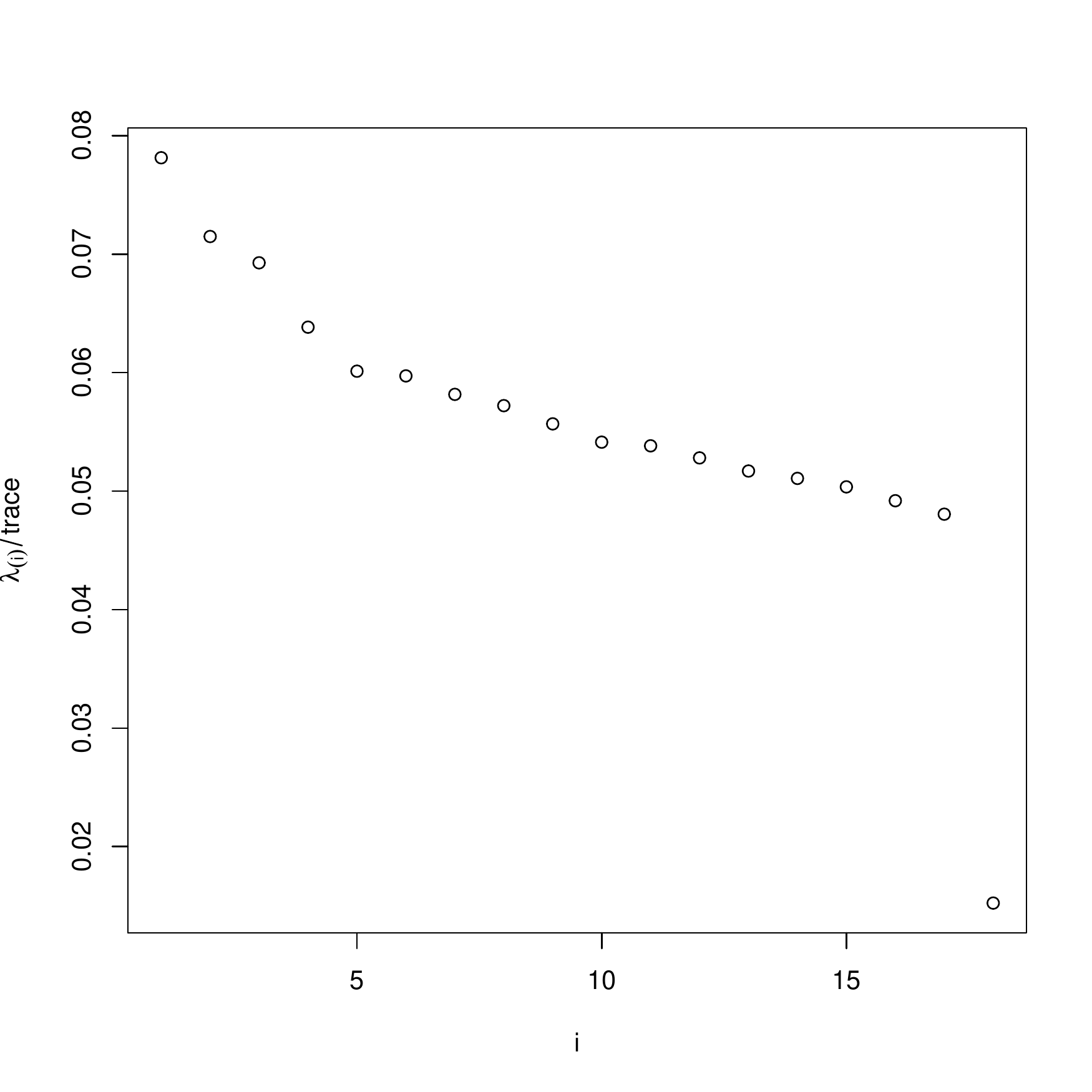}
\includegraphics[width=6.5cm]{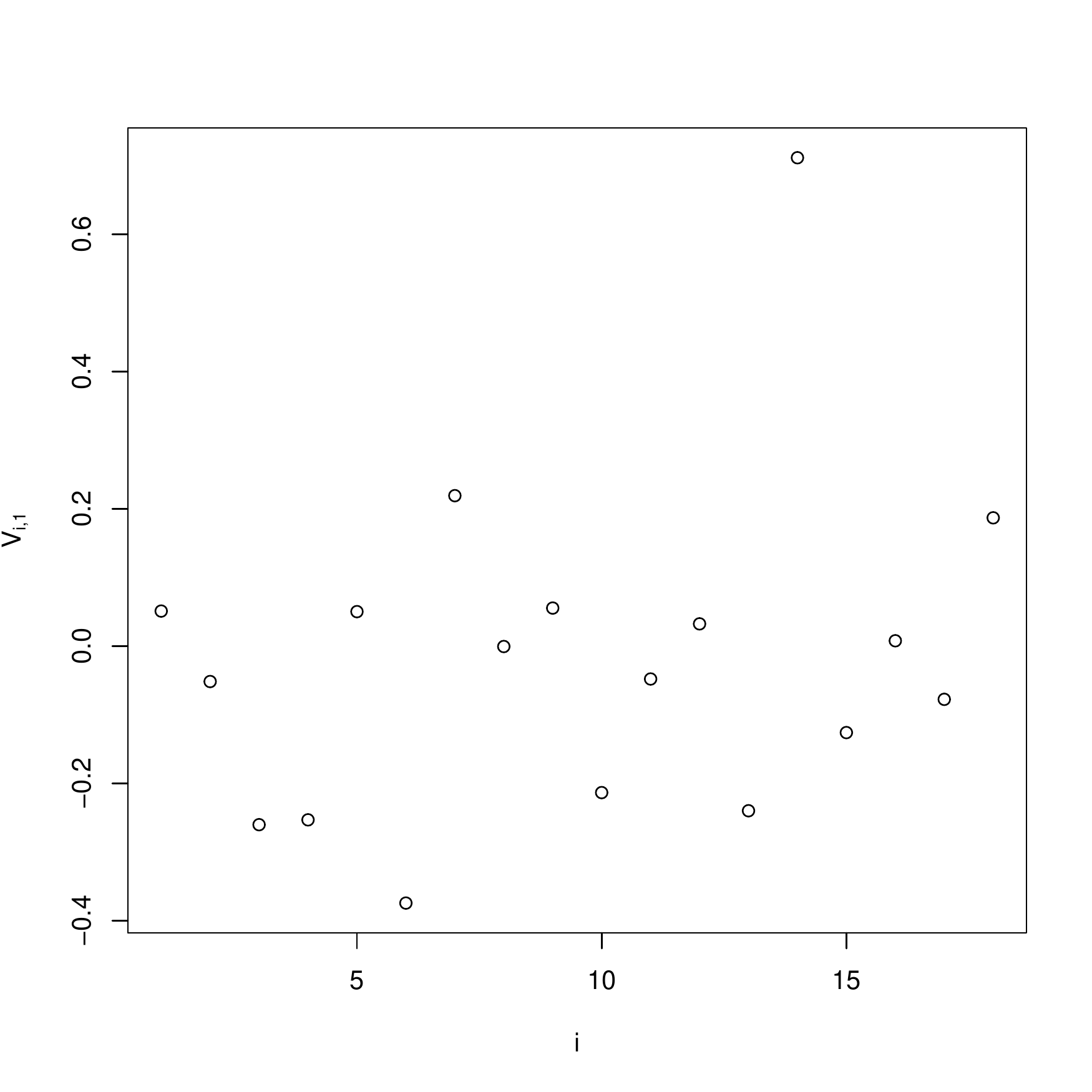}}
\end{figure}

Figure~\ref{eigen:data} shows the ordered eigenvalues of the sample covariance matrix (normalized by its trace) 
and the eigenvector of the FX rate data corresponding to the largest eigenvalue.
There exists a notable spectral gap 
between the largest and second largest eigenvalues and the unit eigenvector corresponding to the largest eigenvector 
has all positive and non-vanishing components.
For comparison and to illustrate the variety of the models discussed above we also plot corresponding realizations of three model specifications from Sections~\ref{sec:case1} and \ref{Sec:case2}. 
In all cases we choose $p=18$ and $n=1567$ in accordance with the data example. 
We assume throughout a moving average structure in the log-volatility process $\log \sigma_{it}$ in \eqref{eq:1a}. More specifically,
\begin{equation}\label{Eq:MA18} \sigma_{it}=\exp(\sum_{k=1}^{18}\eta_{i-k,t}), \;\;\; 1 \leq i \leq 18, \; t \in \mathbb{Z}. \end{equation}
In accordance with the model properties discussed in Section \ref{sec:case1}, we first assume iid\ standard Gaussian $\eta_{i,t}$ and iid\ $Z_{it}$ with a Student-$t$ distribution with $t=3$ degrees of freedom. 
Figure \ref{eigen:model1} shows the normalized eigenvalues and the first unit eigenvector from a realization of this model. We notice a relatively large gap between the first and second eigenvalue and, in accordance with Section \ref{Subsub:eigenvectors}, we see that the first unit eigenvector is relatively close to a unit basis vector. 
Figure \ref{eigen:model2} shows the corresponding realizations for the model \eqref{Eq:MA18} with a specification according to Example \ref{exam:1}, i.e.,\ Exponential(3)-distributed iid\ $\eta_{i,t}$ (meaning that $\P(\eta_{i,t}>x)=\exp(-3x), x \geq 0$, 
which implies $\alpha=3$ and $\E[\ex^{3 \eta}]=\infty$) and iid\ standard Gaussian $Z_{it}$. Compared to the first simulated model, we see a slower decay in the magnitude of the ordered eigenvalues and a more spread out first unit eigenvector. This 
observation illustrates that although the limit behavior of this model and the one analyzed 
before should be very similar (cf. Example \ref{exam:1}), convergence to the prescribed limit appears slower for the heavy-tailed volatility sequence than for the heavy-tailed innovations. 
Finally, Figure \ref{eigen:model3} shows a simulation drawn from \eqref{Eq:MA18} where the $\eta_{i,t}$ are iid\ such that $\P(\eta_{i,t}>x) \sim x^{-2}\exp(-3x), x \to \infty,$ and the $Z_{it}$ are iid\ standard Gaussian. Again, $\alpha=3$, but direct calculations show that the distribution of $\eta_{i,t}$ is convolution equivalent, i.e.,\ it satisfies \eqref{Eq:conv:eq:eta} instead of \eqref{Eq:inf:mom:eta}. The graphs are in line with the analysis in Example \ref{exam:2} and illustrate a very spread out dominant eigenvector. We note that while none of the three very simple models analyzed in the simulations above is able to fully describe
the behavior of the analyzed data, the two models with heavy-tailed volatility and light-tailed innovations are able 
to explain a non-concentrated first unit eigenvector of the sample covariance matrix and 
therefore non-negligible dependence between components as seen in the data.

\section*{Acknowledgements}
Thomas Mikosch's and Xiaolei Xie's research is partly supported by the Danish Research Council Grant DFF-4002-00435 ``Large random matrices with heavy tails and dependence''. Parts of the paper were written when Mohsen Rezapour visited the Department of Mathematics at the University of Copenhagen December 2015--January 2016. He would like to thank the Department of Mathematics for hospitality.

\section{Some $\alpha$-stable limit theory}\setcounter{equation}{0}\label{App:A}
In this paper, we make frequently use of
Theorem 4.3 in Mikosch and Wintenberger \cite{mikosch:wintenberger:2016} which we quote for convenience:
\begin{theorem}\label{thm:mikwin}
Let $(\bfY_t)$ be an $\bbr^p$-valued strictly stationary \seq , $\bfS_n=\bfY_1+\cdots +\bfY_n$ and $(a_n)$ be \st\ $n\,\P(\|\bfY\|>a_n)\to 1$.
Also write for $\vep>0$, $\ov \bfY_t = \bfY_t \1(\|{\mathbf Y}_t\|\le \vep a_n)$, $\underline \bfY_t= \bfY_t-\ov  \bfY_t $ and
\beao
\ov \bfS_{l,n}= \sum_{t=1}^l \ov \bfY_t \qquad  \un \bfS_{l,n}= \sum_{t=1}^l \un \bfY_t\,.
\eeao
Assume the following conditions:
\begin{enumerate}
\item 
$(\bfY_t)$ is \regvary\ with index $\alpha\in (0,2)  \setminus\{1\}$ and spectral tail process $(\bfTh_j)$.
\item
A mixing condition holds: there exists an integer \seq\ $m_n\to\infty$ \st\ $k_n= [n/m_n]\to \infty$
and 
\beam\label{eq:chfa}
\E \ex^{i\bft'\underline \bfS_n/a_n} - \Big(\E \ex^{i\bft'\underline \bfS_{m_n,n}/a_n}\Big)^{k_n} \to 0\,,\qquad \nto\,,\qquad \bft\in\bbr^p\,.
\eeam
\item An anti-clustering condition holds:
\beam\label{eq:acl}
\lim_{l\to\infty} \limsup_{\nto} \P\big(\max_{t=l,\ldots,m_n} \|\bfY_t\|>\delta a_n\mid \|\bfY_0\|>\delta a_n\big)=0\,,\qquad \delta>0\,
\eeam
for the same sequence $(m_n)$ as in (2).
\item If $\alpha\in (1,2)$, in addition $\E[\mathbf{Y}]=\mathbf{0}$ and the vanishing small values condition holds: 
\beam\label{eq:vansm}
\lim_{\vep\downarrow 0}\limsup_{\nto} \P\big(a_n^{-1} \|\ov \bfS_n-\E [\ov \bfS_n]\|>\delta \big)=0\,,\qquad \delta>0\,
\eeam
and $\sum_{i=1}^\infty \E[\|\bfTh_i\|]<\infty$.
\end{enumerate}
Then $a_n^{-1} \bfS_n\std \xi_\alpha$ for an $\alpha$-stable $\bbr^p$-valued vector $\xi_\alpha$ with log-\chf
\beam\label{eq:chfid}
\int_0^\infty \E\big[\ex^{i\,y\,\bft'\sum_{j=0}^\infty \bfTh_j}- \ex^{i\,y\,\bft'\sum_{j=1}^\infty \bfTh_j}-i\,y\,\bft'\1_{(1,2)}(\alpha)\big]\,
d(-y^{\alpha})\,,\qquad \bft\in\bbr^p\,.
\eeam
\end{theorem}\noindent 
\begin{remark}\label{rem:case:alpha:1}
If we additionally assume that $\mathbf{Y}$ is symmetric, which implies $\E[\ov \bfY]=\mathbf{0}$, then the statement 
of the theorem also holds for $\alpha=1$.
\end{remark}

\section{(Joint) Tail behavior for products of \regvary\ \\ \rv s}\label{App:B}
In this paper, we make frequently use of the tail behavior of products of non-negative independent \rv s $X$ and $Y$.
In particular, we are interested in conditions for the existence of the limit
\beam\label{eq:rr}
\lim_{\xto}\dfrac{\P(XY>x)}{\P(X>x)}=q\,.
\eeam 
for some $q\in [0,\infty]$.
We quote some of these results
for convenience.
\ble\label{lem:product} Let $X$ and $Y$ be independent \rv s.
\begin{enumerate}
\item 
If $X$ and $Y$ are \regvary\ with index $\alpha>0$ then $XY$ is \regvary\ with the same index.
\item
If $X$ is \regvary\ with index $\alpha>0$ and $\E[Y^{\alpha+\vep}]<\infty$ for some $\vep>0$ then \eqref{eq:rr} holds with 
$q=\E[Y^\alpha]$.
\item
If  $X$ and $Y$ are iid \regvary\ with index $\alpha>0$ and $\E[Y^\alpha]<\infty$, then \eqref{eq:rr} holds with
$q=2 \E[|Y|^\alpha]$ iff
\beam\label{eq:hh}
\lim_{M\to\infty}\limsup_{\xto}\dfrac{\P(XY>x,M<Y\le x/M )}{\P(X>x)}=0\,.
\eeam
\item
If $X$ and $Y$ are \regvary\ with index $\alpha>0$, $\E[Y^\alpha+ X^\alpha]<\infty$, $\lim_{\xto} \P(Y>x)/\P(X>x)=0$ and 
\eqref{eq:hh} holds, then \eqref{eq:rr} holds with $q=\E[|Y|^\alpha]$.
\item
Assume that $\E[|Y|^\alpha]=\infty$. Then  \eqref{eq:rr} holds with $q=\infty$.
\end{enumerate}
\ele
\begin{proof} (1) This is proved in Embrechts and Goldie \cite{embrechts:goldie:1980}.\\
(2) This is Breiman's \cite{breiman:1965} result.\\
(3) This is Proposition~3.1 in Davis and Resnick \cite{davis:resnick:1985}.\\
(4) This part is proved similarly to (3); we borrow the ideas from \cite{davis:resnick:1985}.
For $M>0$ we have the following decomposition
\beao
\dfrac{\P(XY>x)}{\P(X>x)}&=& \dfrac{\P(XY>x,Y\le M)}{\P(X>x)}+\dfrac{\P(XY>x,M<Y\le x/M)}{\P(X>x)}
+\dfrac{\P(XY>x, Y>x/M)}{\P(X>x)}\\
&\sim & \E[Y^\alpha \1(Y\le M)]+ \dfrac{\P(XY>x,M<Y\le x/M)}{\P(X>x)} + \E [(X\wedge M)^\alpha] \dfrac{\P(Y>x)}{\P(X>x)}\\
&=&  \E[Y^\alpha \1(Y\le M)] + \dfrac{\P(XY>x,M<Y\le x/M)}{\P(X>x)} +o(1)\,.
\eeao
Here we applied Breiman's result twice. The second term vanishes by virtue of \eqref{eq:hh}. Thus $q=  \E[Y^\alpha]$.\\
(5) The same argument as for (4) yields as $\xto$,
\beao
\dfrac{\P(XY>x)}{\P(X>x)}&\ge & \dfrac{\P(XY>x,Y\le M)}{\P(X>x)}\sim \E[Y^\alpha \1(Y\le M)]\,.
\eeao  
Then \eqref{eq:rr} with $q=\infty$ is immediate.
\end{proof}
\ble\label{Lem:splitup} Let $Y_1, \ldots, Y_p \geq 0$ be iid regularly varying \rv s with index $\alpha>0$. Assume that
\begin{equation}\label{Eq:conv:equivalent} \lim_{t \to \infty} \frac{\P(Y_1\cdot Y_2>t)}{\P(Y_1>t)}=c \in (0,\infty)\,.
\end{equation}
Then for any $a_1, \ldots, a_p \geq 0$ such that $a_{\max}:=\max_{j=1, \ldots, p}a_j>0$ and any $v>0$ we have
\begin{eqnarray}\label{Eq:lemma:prod:1} \lim_{t \to \infty} \frac{\P(\prod_{i=1}^p Y_i^{a_i}>vt)}{\P(Y_1^{a_{\max}}>t)}= \sum_{j:a_j=a_{\max}} \lim_{s \to 0} \lim_{t \to \infty} \frac{\P(\prod_{i=1}^p Y_i^{a_i}>vt, Y_j^{a_{\max}}>st)}{\P(Y_1^{a_{\max}}>t)} \end{eqnarray}
and
\begin{eqnarray}\label{Eq:lemma:prod:2}  \lim_{s \to 0} \limsup_{t \to \infty}\frac{\P(\prod_{i=1}^pY_i^{a_i}>vt, \max_{j=1, \ldots, p}Y_j^{a_{\max}}\leq st) }{\P(Y_1^{a_{\max}}>t)}=0. \end{eqnarray}
\begin{proof} 
In view of Davis and Resnick \cite{davis:resnick:1986} the only possible value for $c$ in \eqref{Eq:conv:equivalent} is $2 \E[Y_1^\alpha]$ (which implies that $\E[Y_1^\alpha]<\infty$). Furthermore, we note that the product $\prod_{j:a_j=a_{\max}}Y_j^{a_j}$ is regularly varying with index $-\alpha/a_{\max}$; see Embrechts and Goldie \cite{embrechts:goldie:1980}, Corollary on p.~245. By Breiman's lemma this implies that
\begin{eqnarray*}
&& \lim_{t \to \infty} \frac{\P(\prod_{i=1}^p Y_i^{a_i}>vt)}{\P(Y_1^{a_{\max}}>t)}\\
&=& \lim_{t \to \infty} \frac{\P(Y_1^{a_{\max}}>vt)}{\P(Y_1^{a_{\max}}>t)}\frac{\P(\prod_{i=1}^p Y_i^{a_i}>vt)}{\P(Y_1^{a_{\max}}>vt)}\\
&=& v^{-\alpha/a_{\max}}  \Big( \prod_{j:a_j \neq a_{\max}} \E[Y_j^{\alpha a_j/a_{\max}}] \Big)\lim_{t \to \infty} \frac{\P(\prod_{j:a_j=a_{\max}} Y_j^{a_{\max}}>vt)}{\P(Y_1^{a_{\max}}>vt)}.
\end{eqnarray*}
By Lemma 2.5 in Embrechts and Goldie \cite{embrechts:goldie:1982} 
(cf.\ also  Chover, Ney and Wainger \cite{chover:ney:wainger:1973}) this equals 
$$ v^{-\alpha/a_{max}} \Big(\prod_{j:a_j\neq a_{\max}} \E[Y_j^{\alpha a_j/a_{\max}}] \Big) |\{j:a_j = a_{\max} \}|\,\E[Y_1^\alpha]^{ |\{j:a_j = a_{\max} \}|-1}. $$
On the other hand, we have
\begin{eqnarray*}
&& \sum_{j:a_j=a_{\max}} \lim_{s \to 0} \lim_{t \to \infty} \frac{\P(\prod_{i=1}^p Y_i^{a_i}>vt, Y_j^{a_{\max}}>st)}{\P(Y_1^{a_{\max}}>t)} \\
&=& \sum_{j:a_j=a_{\max}} \lim_{s \to 0} \lim_{t \to \infty} \frac{\P(Y_j^{a_{\max}}\min(s^{-1},v^{-1}\prod_{k\neq j}Y_k^{a_k})>t)}{\P(Y_j^{a_{\max}}>t)} \\
&=&  \sum_{j:a_j=a_{\max}} \lim_{s \to 0} \E[(\min(s^{-1},v^{-1}\prod_{k \neq j}Y_k^{a_k}))^{\alpha/a_{\max}}] \\
&=&  v^{-\alpha/a_{\max}} \sum_{j:a_j=a_{\max}} \prod_{k \neq j}\E[Y_k^{\alpha a_k/a_{\max}}] \\
&=& v^{-\alpha/a_{\max}}  \Big(\prod_{j:a_j\neq a_{\max}} \E[Y_j^{\alpha a_j/a_{\max}}] \Big)  |\{j: a_j=a_{\max}\}|\, \E[Y_1^\alpha]^{|\{j:a_j = a_{\max} \}|-1},
\end{eqnarray*}
where we applied Breiman's lemma in the second step to the bounded random variable $\min(s^{-1},$ $v^{-1}\prod_{k\neq j}Y_k^{a_k})$, and the monotone convergence theorem in the penultimate step. This proves \eqref{Eq:lemma:prod:1}. To prove \eqref{Eq:lemma:prod:2} note that 
for $s>0$,
\begin{eqnarray*}  \frac{\P(\prod_{i=1}^pY_i^{a_i}>vt) }{\P(Y_1^{a_{\max}}>t)}  &\geq& \frac{\P(\prod_{i=1}^pY_i^{a_i}>vt, \max_{j=1, \ldots, p}Y_j^{a_{\max}}\leq st) }{\P(Y_1^{a_{\max}}>t)} \\
&& + \sum_{j:a_j =a_{\max}} \frac{\P(\prod_{i=1}^pY_i^{a_i}>vt, Y_j^{a_{\max}}>st) }{\P(Y_1^{a_{\max}}>t)}\\
&& -\frac{\P(\prod_{i=1}^pY_i^{a_i}>vt, Y_{j_1}^{a_{\max}}>st, Y_{j_2}^{a_{\max}}>st \mbox{ for some } j_1\neq j_2) }{\P(Y_1^{a_{\max}}>t)}, \;\; s>0.
\end{eqnarray*} 
The last summand on the \rhs\  
converges to 0 as $t \to \infty$ by independence of the $Y_j's$. Moreover,  the left-hand term and the second term 
on the \rhs\ become equal by first $t \to \infty$ and then $s \to 0$, in view of \eqref{Eq:lemma:prod:1}. Therefore 
the first right-hand term vanishes by first $t \to \infty$ and then $s \to 0$. This proves the statement.
\end{proof}
\ele
\bpr\label{Pr:genRVforproducts} Let $Y_1, \ldots, Y_p \geq 0$ be iid\ regularly varying with index $\alpha$ and 
$(a_{ij}) \in [0,\infty)^{n\times p}, n,p \geq 1,$ be such that $\max_{1\leq i \leq n} a_{ik}=a_{\max}:=\max_{i,j}a_{ij}>0$ for any $1 \leq k \leq p$. 
\begin{itemize}
\item[(i)] Assume that \eqref{Eq:conv:equivalent} holds. Then the random vector
\begin{equation}\label{Eq:def:vector:Y} \mathbf{Y}:=\big(\prod_{j=1}^pY_j^{a_{ij}}\big)_{1 \leq i \leq n}
\end{equation}
is regularly varying with index $\alpha/a_{\max}$. Furthermore, up to a constant the limit measure $\mu$ 
of $\bfY$ is given by  $\sum_{j=1}^p \mu_j,$ where for any Borel set $B \in [0,\infty]^n$ bounded away from $\mathbf{0}$ and $\nu_\alpha(dz)= \alpha z^{-\alpha-1}dz$,
\begin{eqnarray}
\label{Eq:RV:prod:1}\mu_j(B)&=&\int_0^\infty \P\bigg(\bigg(\1(a_{ij}=a_{\max})z^{a_{\max}}\prod_{k \neq j}Y_k^{a_{ik}}\bigg)_{1 \leq i \leq n} \in B\bigg)\nu_{\alpha}(dz). 
\end{eqnarray}
\item[(ii)] Assume that $\E[Y_1^\alpha]=\infty$. Set 
\beao
p_{\mbox{\scriptsize eff}}&:=&\max_i |\{1\leq j \leq p:a_{ij}=a_{\max}\}|\,,\\
P_{\mbox{\scriptsize eff}}&:=&\{A \subset \{1, \ldots, p\}: |A|=p_{\mbox{\scriptsize eff}} \wedge \, \exists \, i:\, \forall\, j \in A: a_{ij}=a_{\max}\}\,.
\eeao Then the random vector $\mathbf{Y}$ in \eqref{Eq:def:vector:Y}  is regularly varying with index $\alpha/a_{\max}$. 
Furthermore,  up to a constant the limit measure $\mu$ of $\bfY$ is equal to $\sum_{A \in P_{\mbox{\scriptsize eff}}} \mu_A,$ where 
for any Borel set $B \in [0,\infty]^n$ bounded away from $\mathbf{0}$,
\begin{eqnarray}
\label{Eq:RV:prod:2} \mu_A(B)&=&\int_0^\infty \P\bigg(\bigg(\1(a_{ij}=a_{\max} \, \forall\, j \in A)z^{a_{\max}}\prod_{k \notin A}Y_k^{a_{ik}}\bigg)_{1 \leq i \leq n} \in B\bigg)\nu_{\alpha}(dz)\,.
\end{eqnarray}

\end{itemize}
\epr
\begin{proof}% \leavevmode
%\begin{itemize}
(i)~Let $B \in [0,\infty]^n$ be a Borel set bounded away from $\mathbf{0}$. For $s>0$ we have
               \begin{eqnarray}
                \label{Eq:three:summands}  \frac{\P(\mathbf{Y} \in tB)}{\P(Y_1^{a_{\max}}>t)} &=&  \frac{\P(\mathbf{Y} \in tB, \max_{j=1, \ldots, p}Y_j^{a_{\max}}\leq st)}{\P(Y_1^{a_{\max}}>t)}
                 + \sum_{j=1}^p  \frac{\P(\mathbf{Y} \in tB, Y_j^{a_{\max}}>st)}{\P(Y_1^{a_{\max}}>t)} \nonumber\\
                && - \frac{\P(\mathbf{Y} \in tB, Y_{j_1}^{a_{\max}}>st, Y_{j_2}^{a_{\max}}>st, \mbox{ for some } j_1 \neq j_2)}{\P(Y_1^{a_{\max}}>t)}.
               \end{eqnarray}
Since $B$ is bounded away from $\mathbf{0}$, there exists $v>0$ and $1 \leq i \leq n$ such that $B \subset \{(x_1, \ldots, x_n) \in [0,\infty]^n: x_i>v\}$. From Lemma \ref{Lem:splitup}, \eqref{Eq:lemma:prod:2} the first summand in \eqref{Eq:three:summands} therefore tends to 0 by first $t \to \infty$ and then $s \to 0$. Furthermore, the third summand converges to zero as $t \to \infty$ by independence of the $Y_j's$. We are thus left to show
$$ \lim_{s \searrow 0}\lim_{t \to \infty} \frac{\P(\mathbf{Y} \in tB, Y_j^{a_{\max}}>st)}{\P(Y_1^{a_{\max}}>t)}=\mu_j(B), \;\;\; 1 \leq j \leq p, $$
with $\mu_j$ as in \eqref{Eq:RV:prod:1}. For $s>0$ write
\begin{eqnarray*} && \lim_{t \to \infty}\frac{\P(\mathbf{Y} \in tB, Y_j^{a_{\max}}>st)}{\P(Y_1^{a_{\max}}>t)}\\
&=& s^{-\alpha/{a_{\max}}} \lim_{t\to \infty} \P(\mathbf{Y} \in tB \mid Y_j^{a_{\max}}>st) \\
&=& s^{-\alpha/{a_{\max}}} \lim_{t \to \infty} \P\left(\left(\left(\frac{Y_j^{a_{\max}}}{st}\right)^{\frac{a_{ij}}{a_{\max}}}s^{\frac{a_{ij}}{a_{\max}}}t^{\frac{a_{ij}}{a_{\max}}-1}\prod_{k\neq j} Y_k^{a_{ik}}\right)_{1\leq i \leq n} \in B \; \bigg| \; Y_j^{a_{\max}}>st \right)\\
&=& s^{-\alpha/{a_{\max}}} \int_1^\infty \P\left(\left(\1(a_{ij}=a_{\max})sy\prod_{k \neq j}Y_k^{a_{ik}}\right)_{1 \leq i \leq n} \in B\right)\nu_{\alpha/a_{\max}}(dy).
\end{eqnarray*}
Substituting $sy$ by $z$ in the integral finally gives
$$ \lim_{s \searrow 0}\lim_{t \to \infty} \frac{\P(\mathbf{Y} \in tB, Y_j^{a_{\max}}>st)}{\P(Y_1^{a_{\max}}>t)}=\int_0^\infty \P\left(\left(\1(a_{ij}=a_{\max})z^{a_{\max}}\prod_{k \neq j}Y_k^{a_{ik}}\right)_{1 \leq i \leq n} \in B\right)\nu_{\alpha}(dz). $$
(ii) Note first that under our assumptions for any $1\leq n_1 < n_2 \leq p$,
\begin{eqnarray}\nonumber \lim_{t \to \infty}\frac{\P(\prod_{j=1}^{n_2} Y_j>t)}{\P(\prod_{j=1}^{n_1} Y_j>t)}
\label{Eq:non:conv:equiv} &=& \lim_{t \to \infty} \int_0^\infty \frac{\P(\prod_{j=1}^{n_1} Y_j>t/y)}{\P(\prod_{j=1}^{n_1} Y_j>t)}P^{\prod_{j=n_1+1}^{n_2}Y_j}(dy)\\
&\geq& \E\left[\prod_{j=n_1+1}^{n_2}Y_j^\alpha\right]=\infty
\end{eqnarray}
by Fatou's lemma and the regular variation of $\prod_{j=1}^{n_1}Y_j$. Write now
\begin{equation}\label{Eq:split:sum}\mathbf{Y}=\sum_{\substack{1 \leq i \leq n \\ |\{j:a_{ij}=a_{\max}\}|=p_{\mbox{\scriptsize eff}}}} \prod_{j=1}^p Y_j^{a_{ij}} \mathbf{e}_i + \sum_{\substack{1 \leq i \leq n \\ |\{j:a_{ij}=a_{\max}\}|<p_{\mbox{\scriptsize eff}}}} \prod_{j=1}^p Y_j^{a_{ij}} \mathbf{e}_i,
\end{equation}
where $\mathbf{e}_i$ stands for the $i$-th unit vector. The first sum can also be written as
\begin{equation}\label{Eq:split:YA}\sum_{A \in P_{\mbox{\scriptsize eff}}}\mbox{diag}((\1(a_{ij}=a_{\max} \, \forall\, j \in A)\prod_{k \notin A}Y_k^{a_{ik}})_{1 \leq i \leq n})\prod_{j \in A}Y_j^{a_{\max}}=:\sum_{A \in P_{\mbox{\scriptsize eff}}}\mathbf{Y}^A, 
\end{equation}
where for each summand the random matrix and the random factor are independent and for the non-zero entries of the matrix we have $a_{ik}<a_{\max}$ since $k \notin A$. Thus, by the multivariate version of Breiman's lemma each $\mathbf{Y}^A$ is a multivariate regularly varying vector with limit measure $\mu_A$ (up to a constant multiplier) as in \eqref{Eq:RV:prod:2} and normalizing function $P(\prod_{i=1}^{p_{\mbox{\scriptsize eff}}}Y_i^{a_{\max}}>x)$. Furthermore, for $A, A' \in P_{\mbox{\scriptsize eff}}$ with $A \neq A'$ and $i,i'$ such that $a_{ij}=a_{\max} \, \forall\, j \in A$ and $a_{i'j}=a_{\max} \, \forall\, j \in A'$ we have
\begin{eqnarray}\label{Eq:asymp:ind:YA} && \frac{\P(\mathbf{Y}^A_i>x,\mathbf{Y}^{A'}_{i'}>x)}{\P(\prod_{i=1}^{p_{\mbox{\scriptsize eff}}}Y_i^{a_{\max}}>x)}\\
\nonumber &=& \frac{\P((\prod_{j \in A \cap A'}Y_j)^{a_{\max}}\prod_{j \in (A \cap A')^c}Y_j^{a_{ij}}>x,(\prod_{j \in A \cap A'}Y_j)^{a_{\max}}\prod_{j \in (A \cap A')^c}Y_j^{a_{i'j}}>x)}{\P(\prod_{i=1}^{p_{\mbox{\scriptsize eff}}}Y_i^{a_{\max}}>x)}\,. 
\end{eqnarray}
By Jan\ss en and Drees \cite{janssen:drees:2016}, Theorem 4.2 (in connection with Remark 4.3 (ii) and the minor 
change that our random variables are regularly varying with index $\alpha$ instead of $1$), the numerator behaves asymptotically like $\P((\prod_{j \in A \cap A'}Y_j)^{a_{\max}}>x)$, since $\kappa_0=a_{\max}^{-1}, \kappa_{j}=0, j \in (A \cap A')^c$ is the unique non-negative optimal solution to
$$ \kappa_0+\sum_{j \in (A \cap A')^c} \kappa_j \to \min ! $$
under
$$\kappa_0 a_{\max}+\sum_{j \in (A \cap A')^c}\kappa_j a_{ij}\geq 1, \;\;\;  \kappa_0 a_{\max}+\sum_{j \in (A \cap A')^c}\kappa_j a_{i'j}\geq 1.$$
This is because $\min(a_{ij},a_{i'j})<a_{\max}$ and $\max(a_{ij},a_{i'j})\leq a_{\max}$ for all $j \in (A \cap A')^c$. Since $A \neq A'$, we have $|A \cap A'|<p_{\mbox{\scriptsize eff}}$ and thus, by \eqref{Eq:non:conv:equiv}, the expression \eqref{Eq:asymp:ind:YA} converges to 0 as $x \to \infty$. Therefore, each component of $\mathbf{Y}^A$ is asymptotically independent of each component of $\mathbf{Y}^{A'}$ and thus the sum in \eqref{Eq:split:YA} is multivariate regularly varying with limit measure $\sum_{A \in P_{\mbox{\scriptsize eff}}}\mu_A$ and normalizing function $\P(\prod_{i=1}^{p_{\mbox{\scriptsize eff}}}Y_i^{a_{\max}}>x)$. Since the second sum in \eqref{Eq:split:sum} consists by \eqref{Eq:non:conv:equiv} only of random vectors for which $\P(\|\prod_{j=1}^p Y_j^{a_{ij}} \mathbf{e}_i\|>x)=\P(\prod_{j=1}^p Y_j^{a_{ij}}>x) =o(\P(\prod_{i=1}^{p_{\mbox{\scriptsize eff}}}Y_i^{a_{\max}}>x))$, we have that $\mathbf{Y}$ is regularly varying with index $\alpha/a_{\max}$ and limit measure $\sum_{A \in P_{\mbox{\scriptsize eff}}}\mu_A$ by Lemma 3.12 in Jessen and Mikosch \cite{jessen:mikosch:2006}. 
\end{proof}

\end{document}